\documentclass[twoside]{article}
\pdfoutput=1
\usepackage[dvipsnames]{xcolor}

\usepackage{graphicx, epstopdf, verbatim, fancyhdr, 
wrapfig}

\usepackage{amsmath}
\input amssym.def
\input amssym.tex
\usepackage{bbm}

\skip\footins 22pt plus 2pt minus 2pt
\newcommand{\ints}{{\mathbbm{Z}}}
\newcommand{\reals}{{\mathbbm{R}}}

\newcommand{\diag}{\mathop{\rm diag}\nolimits}
\newcommand{\sgn}{\mathop{\rm sgn}\nolimits}

\newcommand{\dast}{{\displaystyle{\ast}}}

\renewcommand{\a}{{\bf a}}

\newcommand{\x}{{\bf x}}
\newcommand{\y}{{\bf y}}
\renewcommand{\v}{{\bf v}}
\newcommand{\w}{{\bf w}}
\newcommand{\p}{{\bf p}}
\newcommand{\q}{{\bf q}}

\begin{document}
\thispagestyle{empty}
\large

\bigskip

\centerline{\LARGE\bf Making matrices better:}
\centerline{\bf Geometry and topology of polar and singular value decomposition\rm}

\bigskip
\centerline{\sc Dennis DeTurck, Amora Elsaify, Herman Gluck, Benjamin Grossmann}
\centerline{\sc Joseph Hoisington, Anusha M.\,Krishnan, Jianru Zhang}

\bigskip

\centerline{\bf Abstract\rm}

{\normalsize{Our goal here is to see the space of matrices of a given size from a geometric 
and topological perspective, with emphasis on the families of various ranks
and how they fit together. We pay special attention to the nearest orthogonal
neighbor and nearest singular neighbor of a given matrix, both of which play
central roles in matrix decompositions, and then against this visual backdrop examine
the polar and singular value decompositions and some of their applications.\\
\textbf{MSC Primary}: 15-02, 15A18, 15A23, 15B10; 
\textbf{Secondary}: 53A07, 55-02, 57-02, 57N12, 91B24, 91G30, 92C55.}}

%

\addtolength{\baselineskip}{1pt}

Figure 1 is the kind of picture we have in mind, in which we focus on $3 \times 3$ matrices,
view them as points in Euclidean 9-space $\reals^9$, ignore the zero matrix at the origin,
and scale the rest to lie on the round 8-sphere $S^8(\sqrt{3})$ of radius $\sqrt{3}$,  so as to include
the orthogonal group $O(3)$.\\
\begin{wrapfigure}[15]{r}{0.5\textwidth}
\vspace{-20pt}
\center{\includegraphics[width=0.48\textwidth]{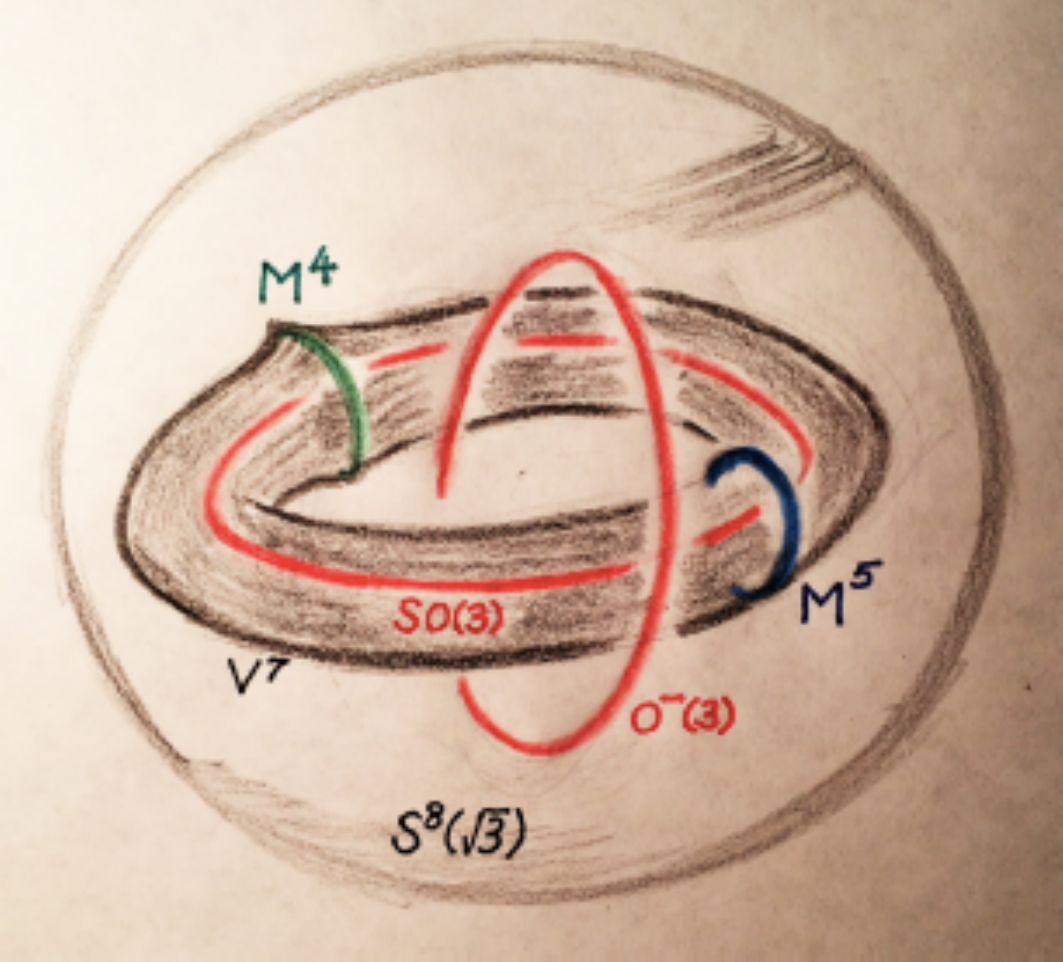}
\caption{{\textit{\textbf{\boldmath A view of $3\times 3$ matrices}}}}
}
\end{wrapfigure}
\noindent The two components of $O(3)$  appear as real projective 3-spaces in the\break 8-sphere, each
the core of a open neighborhood of nonsingular matrices, whose cross-sectional fibres are triangular 
5-dimensional cells lying on great\break 5-spheres.  The common boundary 
of these two neighborhoods is the 7-dimensional algebraic variety $V^7$ of singular matrices.

This variety fails to be a submanifold precisely along the 4-manifold $M^4$ of 
matrices of rank 1. The complement $V^7 -M^4$, consisting of matrices of rank 2, is a large tubular neighborhood of 
a core 5-manifold $M^5$ consisting of the ``best matrices of rank 2'', namely those which are orthogonal on a 
2-plane through the origin and zero on its orthogonal complement.  $V^7$ is filled by geodesics, 
each an eighth of a great circle on the 8-sphere, which run between points of $M^5$ and  $M^4$  
with no overlap along their interiors.  A circle's worth of these geodesics originate from each point of $M^5$,  
leaving it orthogonally, and a 2-torus's worth of these geodesics arrive at each point of $M^4$, also orthogonally.

We will confirm the above remarks, determine the topology and geometry of all these pieces, and the 
interesting cycles (families of matrices) which generate some of their homology, see how they all fit together to form the 
8-sphere, and then in this setting visualize the polar and singular value decompositions and some of their applications.

In Figure 2, we start with a $3\times 3$ matrix $A$ with positive determinant on $S^8(\sqrt{3})$, and show its polar
and singular value decompositions, its nearest orthogonal neighbor $U$, and its nearest singular neighbor
$B$ on that 8-sphere.
\begin{figure}[h!]
\center{\includegraphics[height=180pt]{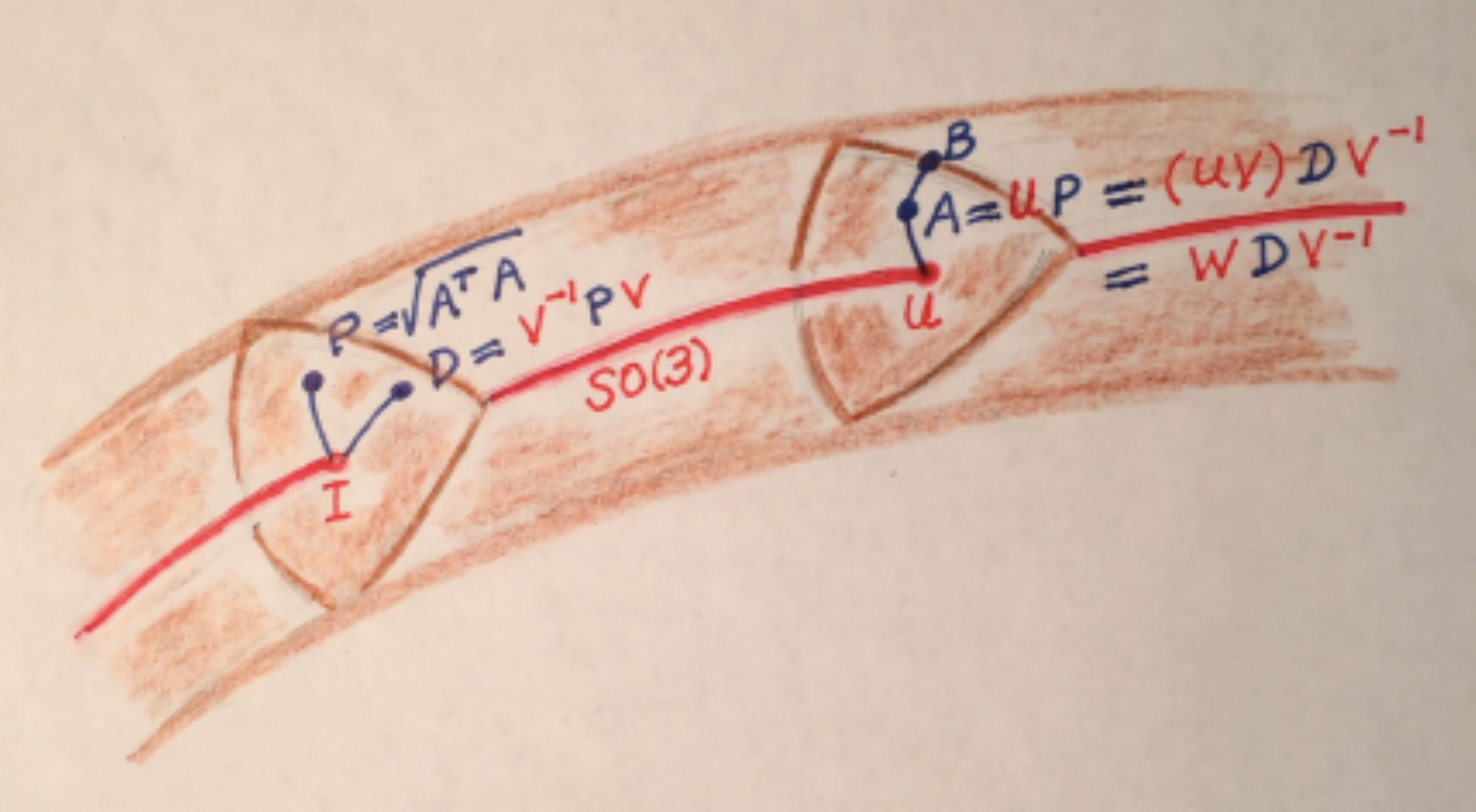}
\caption{{\textit{\textbf{\boldmath Polar and singular value decomposition of $A$}}}}
}
\end{figure}

\noindent Since $\det A>0$, $A$ lies inside the tubular neighborhood $N$ of $SO(3)$ on the 8-sphere.
The nearest orthogonal neighbor $U$ to $A$ is at the center of the 5-cell fibre of $N$ containing $A$,  
while the nearest singular neighbor $B$ to $A$ lies on the boundary of that 5-cell.

\textit{\textbf{These two nearest neighbors play a central role in the applications.}}

The positive definite symmetric matrix  $P  = \sqrt{A^TA} =  U^{-1}A$   lies on the corresponding fibre of  
$N$  centered at the identity $I$.
Orthogonal diagonalization of  $P$  yields the diagonal matrix  $D  =  V^{-1} P V$  on that same fibre, with  
$V  \in  SO(3)$.  

Then we have the two matrix decompositions
\begin{align*}
A&={\color{red}U}\,{\color{blue}P}\qquad\quad\mbox{(polar decomposition)}\\
&={\color{red}U}({\color{red}V}{\color{blue}D}{\color{red}V^{-1}}) 
={\color{red}UV}{\color{blue}D}{\color{red}V^{-1}}
={\color{red}W}{\color{blue}D}{\color{red}V^{-1}} \quad\mbox{(singular value decomposition)}
\end{align*}

Polar and singular value decompositions have a wealth of applications, from which we
sample the following: least squares estimate of satellite attitude as well as computational
comparative anatomy (both instances of nearest orthogonal neighbor, and known as the
\textit{\textbf{Orthogonal Procrustes Problem}}); and facial recognition via eigenfaces as well as interest
rate term structures for US treasury bonds (both instances of nearest singular neighbor and known as
\textit{\textbf{Principal Component Analysis}}).

\bigskip

\noindent\textbf{To the reader}.

In the first half of this paper, we focus on the geometry and topology of spaces
of matrices, quickly warm up with the simple geometry of $2 \times 2$ matrices, and then concentrate entirely
on the surprisingly rich and beautiful geometry of $3 \times 3$ matrices.
Hoping to have set the stage well in that case,
we go no further on to higher dimensions, but invite the inspired reader to do so.

In the second half of the paper, we consider matrices of arbitrary size and shape,
as we focus on their singular value and polar decompositions, and applications of these,
and suggest a number of references for further reading.

As usual, figures depicting higher-dimensional phenomena are at best artful lies, emphasizing
some features and distorting others, and need to be viewed charitably and cooperatively
by the reader.

\bigskip

\noindent\textbf{Acknowledgments}.

We are grateful to our friends Christopher Catone, Joanne Darken, Ellen Gasparovic,
Chris Hays, Kostis Karatapanis, Rob Kusner and Jerry Porter for their help with this paper.

\bigskip

\vfill
\eject

\centerline{\LARGE{\textbf{\textit{Geometry and topology of spaces of matrices}}}}

\bigskip

\centerline{{\Large\textbf{\boldmath{$2\times 2$ matrices}}}}

{\normalsize{We begin with $2\times 2$ matrices, view them as points in Euclidean 4-space  $\reals^4$,
ignore the zero matrix at the origin, and scale the rest to lie on the round 3-sphere
$S^3(\sqrt{2})$  of radius $\sqrt{2}$,  so as to include the orthogonal group $O(2)$.}}

\begin{wrapfigure}[15]{r}{0.5\textwidth}
\vspace{-34pt}
\center{\includegraphics[width=0.48\textwidth]{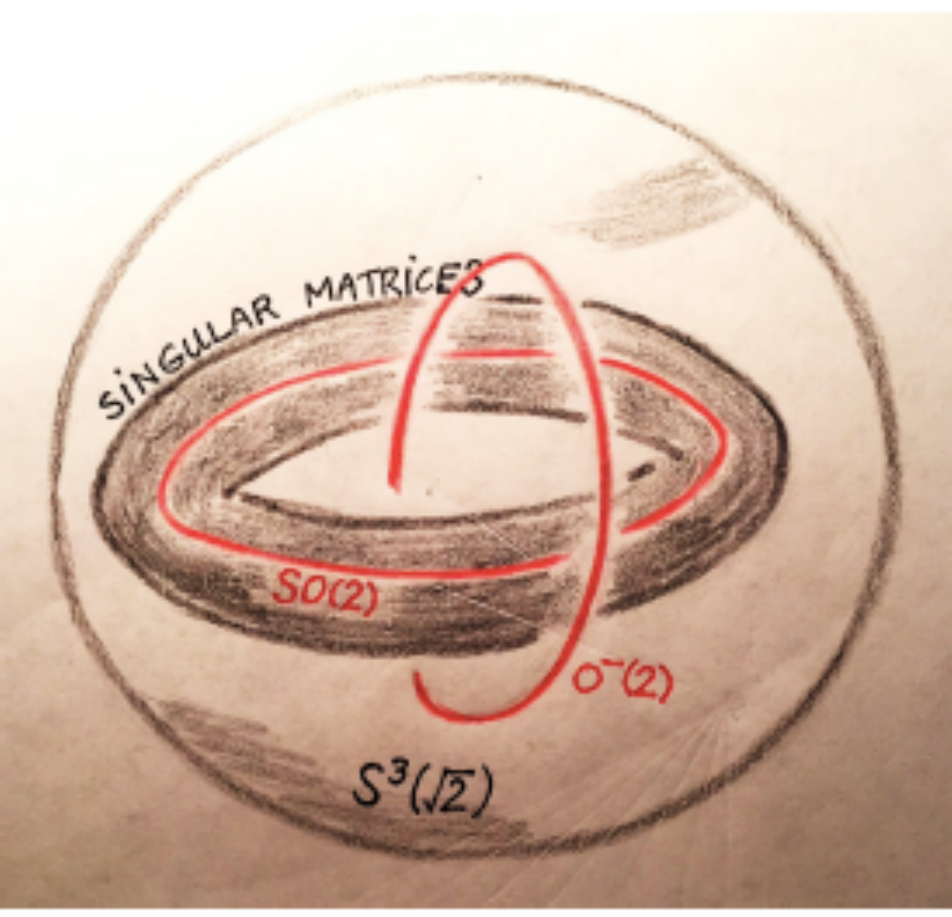}
\caption{{\textit{\textbf{\boldmath A view of $2\times 2$ matrices}}}}
}
\end{wrapfigure}
\noindent\textbf{(1) First view}.  A simple coordinate change reveals that within this 3-sphere, the two components  
$SO(2)$  and  $O^-(2)$  of  $O(2)$  appear  as linked 
orthogonal great circles, while the singular matrices 
appear as the Clifford torus halfway between these two great circles (Figure 3).
The
complement of this Clifford torus consists of open tubular neighborhoods $N$ and $N'$
of $SO(2)$ and $O^-(2)$, each an open solid torus.

\noindent\textbf{(2) Features}.\\
\vspace{-25pt}
\begin{enumerate}
\item[(i)] On  $S^3(\sqrt{2})$,  the determinant function det takes its maximum value of  $+1$
on  $SO(2)$,  its minimum value of  $-1$  on  $O^-(2)$  and its intermediate value of  0
on the Clifford torus of singular matrices.
\end{enumerate}

\vspace{-20pt}
\begin{enumerate}
\item[(ii)] The level sets of  det  on  $S^3(\sqrt{2})$  are tori parallel to the
Clifford torus, and the great circles $SO(2)$   and
$O^-(2)$.
\item[(iii)] The orthogonal trajectories to these level sets (i.e., the gradient flow lines of  det)  are quarter circles 
which leave $SO(2)$  orthogonally and arrive at  $O^-(2)$  orthogonally.  
\item[(iv)] The symmetric matrices on $S^3(\sqrt{2})$ lie on a great 2-sphere with $I$ and $-I$ as poles
and with $O^-(2)$ as equator. Inside the symmetric matrices, the diagonal matrices
appear as a great circle through these poles, passing alternately through the tubular
neighborhoods $N$ and $N'$ of $SO(2)$ and $O^-(2)$, and crossing the Clifford torus four
times.
\item[(v)] On the great 2-sphere of symmetric matrices,
the round disk of angular radius $\pi/4$ centered at  $I$  is one of the cross-sectional fibres of the tubular 
neighborhood  $N$  of  $SO(2)$.  It meets  $SO(2)$ orthogonally at its
center, and meets the Clifford torus orthogonally along its boundary, thanks to (i), (ii) and (iii) above.
\item[(vi)] The tangent space to $S^3(\sqrt{2})$ at the identity matrix $I$ decomposes orthogonally
into the one-dimensional space of skew-symmetric matrices (tangent to $SO(2)$), and the two-dimensional
space of traceless symmetric matrices, tangent to the great 2-sphere of symmetric matrices. Within the traceless symmetric
matrices is the one-dimensional space of traceless diagonal matrices, tangent to the great circle of diagonal 
matrices.
\item[(vii)] Left or right multiplication by elements of  $SO(2)$  are isometries of  $S^3(\sqrt{2})$  
which take this cross-sectional fibre of  $N$  at  $I$  to the corresponding cross-sectional fibres of  $N$  
at the other points along  $SO(2)$.  Left or right multiplication by elements of  
$O^-(2)$  take this fibration of  $N$  to the corresponding fibration of  $N'$.
\end{enumerate}

\medskip

\noindent\textbf{(3)  Nearest orthogonal neighbor}.  Start with a nonsingular $2\times 2$ matrix  $A$  on  $S^3(\sqrt{2})$ 
and suppose, to be specific, that  $A$  lies in the open tubular neighborhood  $N$  
of  $SO(2)$.
\textit{\textbf{\boldmath We 
claim that the nearest orthogonal neighbor to $A$  on that 3-sphere 
is the center of the cross-sectional fibre of  $N$  on which it lies.}}

To see this, note that a geodesic (great circle arc) from  $A$  to its nearest neighbor  $U$
on  $SO(2)$ must meet  $SO(2)$  orthogonally at  $U$,  
and therefore must lie in the cross-sectional fibre of  $N$  through  $U$.  It follows that  $A$  also lies in that fibre, whose center 
is at  $U$,  confirming the above claim.

\begin{figure}[h!]
\center{\includegraphics[height=144pt]{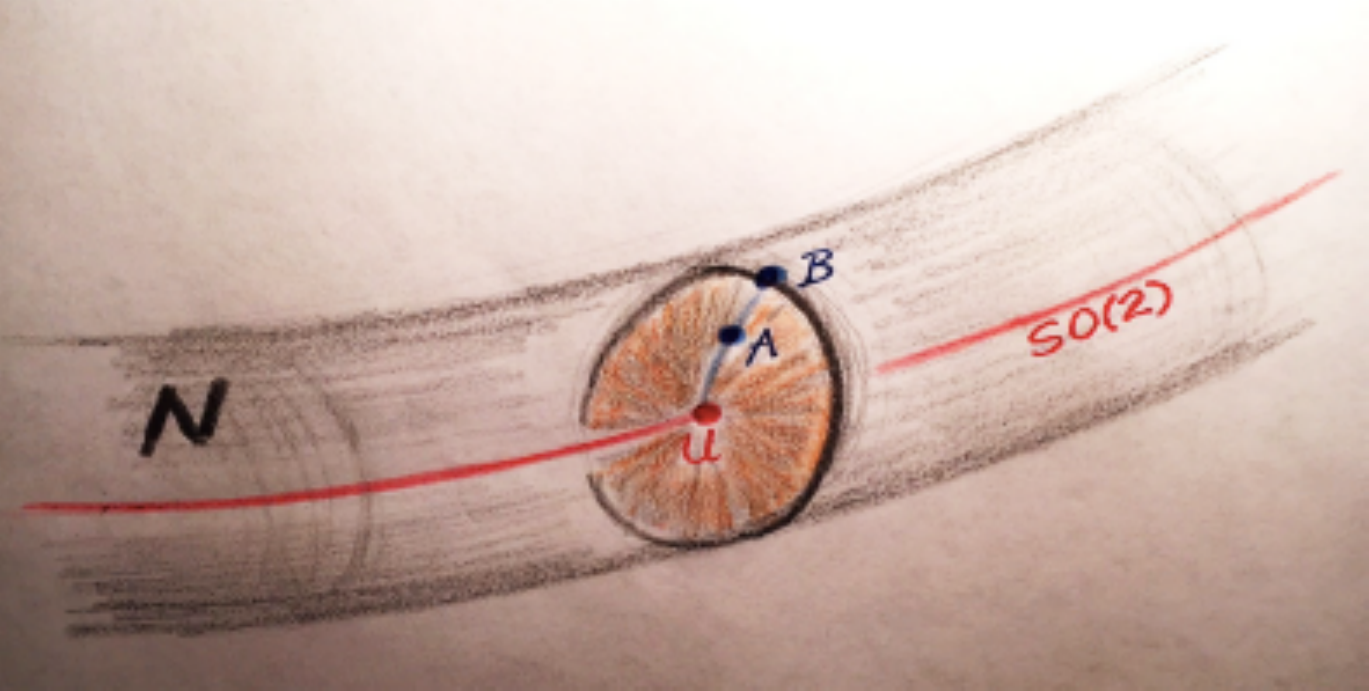}
\caption{{\textit{\textbf{\boldmath Nearest orthogonal and nearest singular neighbors to a matrix $A$}}}}
}
\end{figure}

\medskip

\noindent\textbf{(4)  Nearest singular neighbor}. Start with a nonsingular $2\times 2$ matrix $A$ on  $S^3(\sqrt{2})$. 
\textit{\textbf{\boldmath
We claim that the nearest singular neighbor to $A$  on that $3$-sphere is on 
the boundary of the cross-sectional disk on 
which it lies, at the end of the ray from 
its center through $A$.}}

To see this, recall from (1) that the level surfaces of det on  $S^3(\sqrt{2})$  are tori parallel 
to the Clifford torus, and that their orthogonal trajectories are the quarter 
circles which leave $SO(2)$ orthogonally and arrive 
at  $O^-(2)$  orthogonally.  It follows that the geodesics orthogonal to the Clifford torus lie in the 
cross-sectional disk 
fibres of 
the tubular neighborhoods $N$  and $N'$ of  $SO(2)$  and  $O^-(2)$.  

Now a geodesic (great circle arc) from  $A$  to its nearest singular neighbor  $B$  on the
Clifford torus must meet that torus orthogonally at  $B$,  and hence must lie in one of
these cross-sectional disk fibres (Figure 4).
If  $A$  is not orthogonal, then $B$  lies at the end of the unique ray from the center of 
this fibre through $A$,  and hence is uniquely determined by  $A$.
If  $A$ is orthogonal, then $B$ can lie at the end of any of the rays from the center $A$  
of this fibre, and so every point on the circular boundary of this fibre is a closest
singular neighbor to $A$  on  $S^3(\sqrt{2})$.

\medskip

\noindent\textbf{(5) Gram-Schmidt}.   Having just looked at the geometrically natural map which takes 
a nonsingular $2\times 2$ matrix to its nearest neighbor on the orthogonal group $O(2)$,  
it is irresistable to compare this 
with the Gram-Schmidt orthonormalization procedure. This procedure depends on a choice of 
basis for $\reals^2$,  hence is not ``geometrically natural'', that is to say, not  
$O(2) \times O(2)$  equivariant.  

\begin{wrapfigure}[10]{r}{0.5\textwidth}
\vspace{-34pt}
\center{\includegraphics[width=0.48\textwidth]{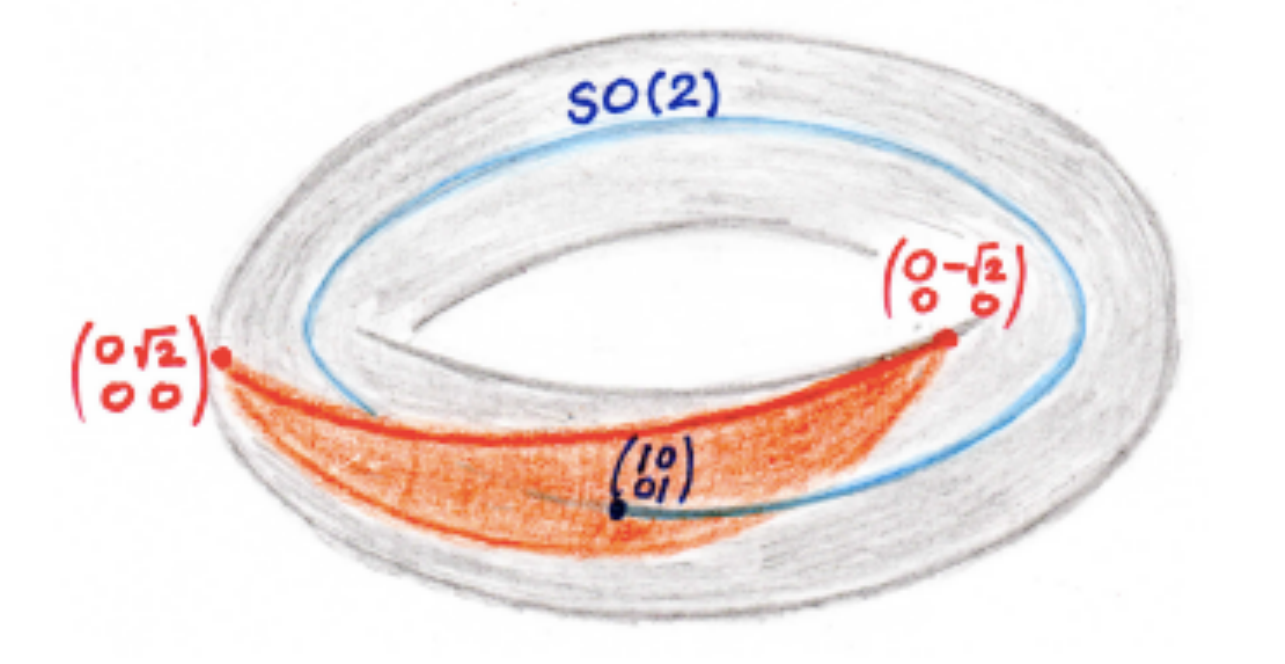}
\caption{{\textit{\textbf{\boldmath $GS^{-1}(I)$ is an open $2$-cell in $S^3(\sqrt{2})$ 
with boundary on the Clifford torus}}}}
}
\end{wrapfigure}
We see this geometric defect in Figure 5, where we restrict the Gram-Schmidt
procedure ${GS}$ to $S^3(\sqrt{2})$, and display the inverse image 
$GS^{-1}(I)$  of the identity  $I$
on that 3-sphere.

The inverse images of the other points on $SO(2)$  are rotated versions of  
$GS^{-1}(I)$.
It is visually evident that this picture, and hence the Gram-Schmidt procedure itself, is not 
equivariant with respect to the 
action of $SO(2)$ via conjugation, which fixes $SO(2)$  pointwise, but rotates $O^-(2)$ within itself.

\bigskip

\bigskip

\centerline{\Large{\textbf{\boldmath{$3\times 3$ matrices}}}}

{\normalsize{We turn now to $3\times 3$ matrices, view them as points in Euclidean 9-space  $\reals^9$,  
once again ignore the zero matrix at the origin, and scale the rest to lie on the 
round 8-sphere $S^8(\sqrt{3})$  of radius $\sqrt{3}$,  so as to include the orthogonal group $O(3)$.}}

\medskip

\noindent\textbf{(1) First view}.  The two components  $SO(3)$  and  $O^-(3)$  of  $O(3)$  appear as real projective 
3-spaces on $S^8(\sqrt{3})$, while the singular matrices (ranks 1 and 2) on this 
8-sphere appear as a 7-dimensional algebraic variety $V^7$ separating them.

\begin{wrapfigure}[14]{r}{0.5\textwidth}
\vspace{-16pt}
\center{\includegraphics[width=0.48\textwidth]{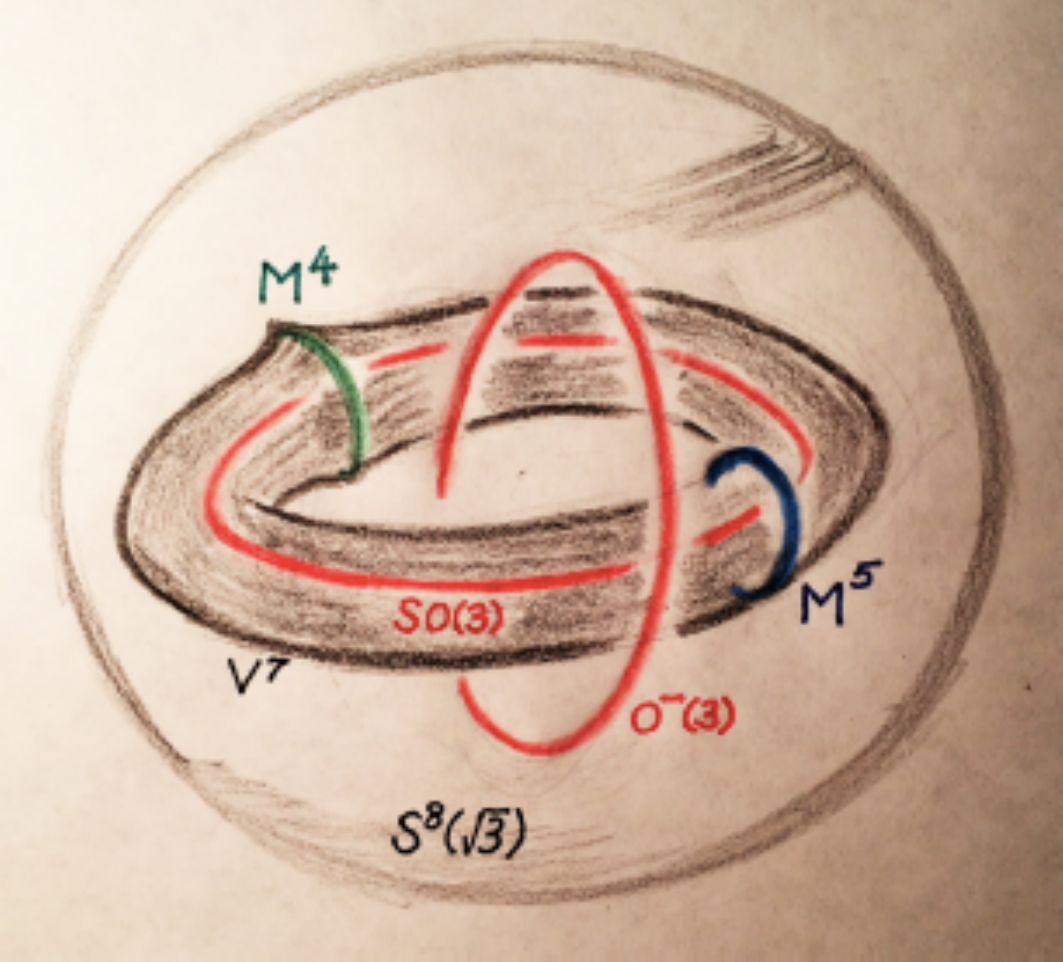}
\caption{{\textit{\textbf{\boldmath A view of $3\times 3$ matrices}}}}
}
\end{wrapfigure}
\noindent Contrary to appearances in Figure 6, the two components of $O(3)$  are too 
low-dimensional to be linked in the 8-sphere.  
The subspaces  $V^7$,  $M^4$  and  $M^5$ in the figure were defined in the introduction, and will be examined
in detail as we proceed.

\medskip

\noindent\textbf{\boldmath (2) The tangent space to $S^8(\sqrt{3})$ at the identity matrix}
decomposes orthogonally
into the three-dimensional space of skew-symmetric matrices (tangent to $SO(3)$), and the five-dimensional
space of traceless symmetric matrices, tangent to the great 5-sphere of symmetric matrices. Within the traceless symmetric
matrices is the two-dimensional space of traceless diagonal matrices, tangent to the great 2-sphere of diagonal 
matrices in $S^8(\sqrt{3})$.

\medskip

\noindent\textbf{\boldmath(3) A 2-sphere's worth of diagonal $3\times 3$ matrices}.
The great 2-sphere of diagonal $3 \times 3$ matrices on  $S^8(\sqrt{3})$  will play a key role in our understanding of the geometry 
of $3\times 3$ matrices as a whole.

\begin{wrapfigure}[16]{r}{0.5\textwidth}
\vspace{-35pt}
\center{\includegraphics[width=0.48\textwidth]{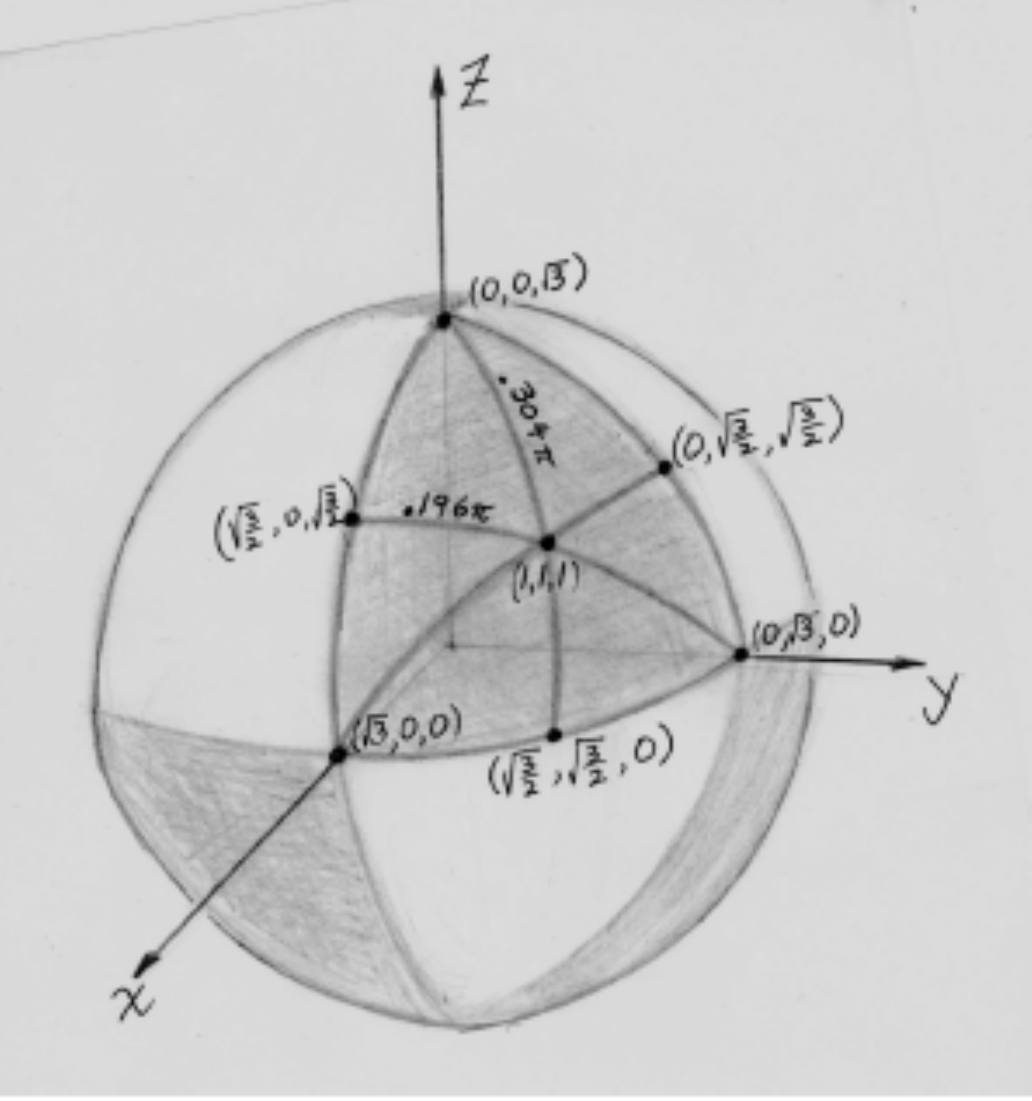}
\caption{{\textit{\textbf{\boldmath Diagonal matrices in $S^8(\sqrt{3})$}}}}
}
\end{wrapfigure}
In Figure 7, the diagonal matrix  $\diag(x, y, z)$  is located at the point  $(x, y, z)$,
and indicated ``distances'' are really angular separations.

This 2-sphere is divided into eight spherical triangles, with the shaded ones centered at the points
$(1, 1, 1)$, $(-1,-1, 1)$, $(1, -1,-1)$ and $(-1, 1,-1)$
of  $SO(3)$,  and the unshaded ones centered at points of  $O^-(3)$.

The interiors of the shaded triangles will lie in the open tubular neighborhood (yet to
be defined) of  $SO(3)$  on  $S^8(\sqrt{3})$,  the interiors of the unshaded triangles will lie in the open tubular neighborhood of 
$O^-(3)$,  while the shared boundaries lie on the variety  $V^7$  
of singular matrices, with the vertices of rank 1, the open edges of rank 2, and the centers of the edges ``best of rank 2''.

\medskip

\noindent\textbf{(4) Symmetries}.  We have  $O(3) \times O(3)$  acting as a group of isometries of our space  $\reals^9$  of all 
$3\times 3$ matrices, and hence of the normalized ones on  $S^8(\sqrt{3})$, via the map
$$(U, V)\,\dast\, A  =  U A V^{-1}.$$
This action is a rigid motion of the 8-sphere which takes the union of the two
$\reals P^3$s  representing $O(3)$ to themselves (possibly interchanging them), and takes 
the variety  $V^7$  of singular matrices separating them to itself.

``Natural geometric constructions'' for $3\times 3$ matrices are those which are equivariant
 with respect to this action of  $O(3)\times O(3)$.

\medskip

\noindent\textbf{(5) \boldmath Tubular neighborhoods of $SO(3)$ and $O^-(3)$ in $S^8(\sqrt{3})$}. We expect, by analogy with 
$2 \times 2$ matrices, that the complement in $S^8(\sqrt{3})$  of the variety $V^7$ of singular matrices consists of open tubular 
neighborhoods of the two components $SO(3)$  and  $O^-(3)$  of the orthogonal group, with fibres which lie 
on the great 5-spheres which meet these cores orthogonally.  

At the same time, our picture of the great 2-sphere's worth of diagonal $3 \times 3$ matrices alerts us that we cannot expect the 
fibres of these neighborhoods to be round 5-cells; instead they must somehow take on the triangular shapes seen in 
Figure 7.

Indeed, look at that figure and focus on the open shaded spherical triangle  $D^2$ 
centered at the identity and lying in the first octant.  Let  $SO(3)$  act on this triangle
by conjugation,
$$A \to U \,\dast\, A  =  U A U^{-1},$$
and the image will be a corresponding open triangular shaped region $D^5$ centered
at the identity on the great 5-sphere of symmetric matrices, and consisting of the
positive definite ones.  Going from  $D^2$  to  $D^5$  is like fluffing up a pillow.

This open 5-cell  $D^5$  is the fibre centered at the identity of the tubular neighborhood  
$N$  of  $SO(3)$,  and the remaining fibres can be obtained by left (say) translation of  
$D^5$  by the elements of $SO(3)$.

Why are these fibres disjoint?  That is, why will two left translates of $D^5$ along
$SO(3)$ be disjoint?

We can see from Figure 7 that it is going to be a close call, 
since the \textit{\textbf{closures}} of the spherical triangles centered at 
$(1, 1, 1)$ and at  $(1, -1, -1)$  meet at the point  $(\sqrt{3}, 0, 0)$,  even though their
interiors are disjoint.

Consider a closed geodesic on  $SO(3)$,  such as the set of transformations  
$$A_t=\left[\begin{array}{ccc}\cos t&-\sin t&0\\\sin t&\cos t&0\\0&0&1\end{array}\right]\qquad 0\le t\le 2\pi.$$

\begin{wrapfigure}[15]{l}{0.45\textwidth}
\vspace{-35pt}
\center{\includegraphics[width=0.43\textwidth]{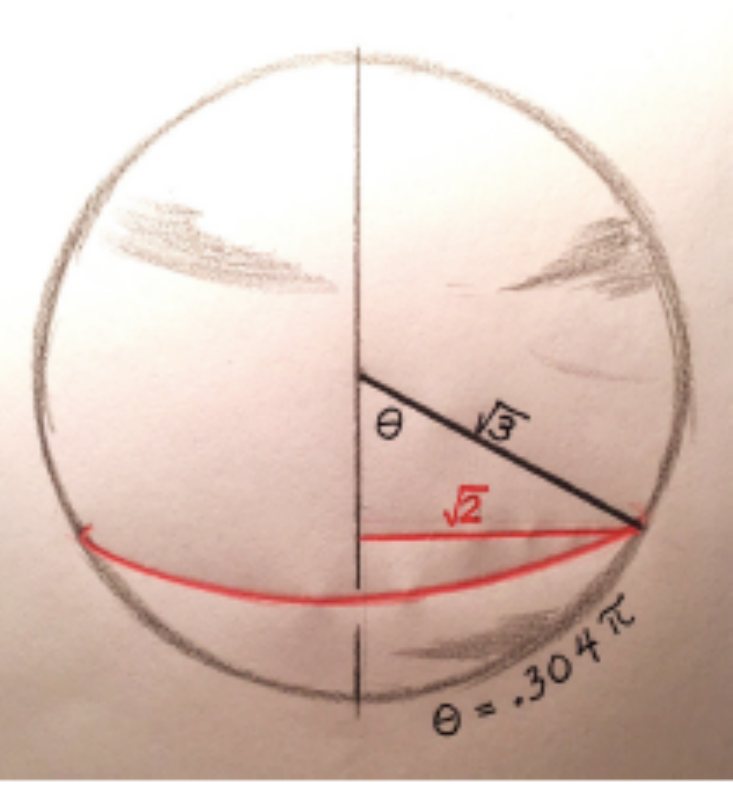}
\caption{{\textit{\textbf{\boldmath A closed geodesic on $SO(3)$ appears as a small circle on $S^8(\sqrt{3})$}}}}
}
\end{wrapfigure}
\noindent
Figure 8 is a picture of that closed geodesic, appearing as a small circle of radius $\sqrt{2}$  on  $S^8(\sqrt{3})$.

In this picture, two great circles which meet the small circle orthogonally will come together at the south pole, after 
traveling an angular distance  $0.304\,\pi$,
but not before.

Since any two points $U$  and  $V$  of  $SO(3)$  lie together on a common closed geodesic (which is a small circle of radius $\sqrt{2}$ on an 
8-sphere of radius $\sqrt{3}$),
and since the maximum angular separation between the center of the 5-disk  $D^5$  and its boundary is  $0.304\,\pi$,  it follows that 
the \textit{\textbf{open}} 5-disks  $UD^5$  and  $VD^5$  must be disjoint.

In this way, we see that the union of the disjoint open 5-disks  $U D^5$,  as  $U$  ranges 
over  $SO(3)$,  forms an open tubular neighborhood  $N$  of  $SO(3)$  in  $S^8(\sqrt{3})$.  This
tubular neighborhood is topologically trivial under the map
$$SO(3) \times D^5 \to N\quad\mbox{via}\quad (U, P) \to U P.$$
In similar fashion, we get an open tubular neighborhood  $N'$  of  $O^-(3)$,  likewise
topologically trivial.
The common boundary of these two tubular neighborhoods is the variety  $V^7$
of singular matrices on  $S^8(\sqrt{3})$.

\medskip

\noindent\textbf{\boldmath (6) The determinant function on $S^8(\sqrt{3})$}.
The determinant function  det  on  $S^8(\sqrt{3})$  takes its maximum value of  $+1$  on  $SO(3)$,  
its minimum value of  $-1$  on  $O^-(3)$,  and its intermediate value of  0  on  $V^7$.

Unlike the situation for $2\times 2$ matrices, the orthogonal trajectories of the level sets
of  det  are not geodesics, since the 5-cell fibres of the tubular neighborhoods  N
and  $N'$  of  $SO(3)$  and  $O^-(3)$  are not round.
In Figure 9, we see the level curves of  det  on the great 2-sphere of diagonal matrices.
\begin{figure}[h!]
\center{\includegraphics[height=180pt]{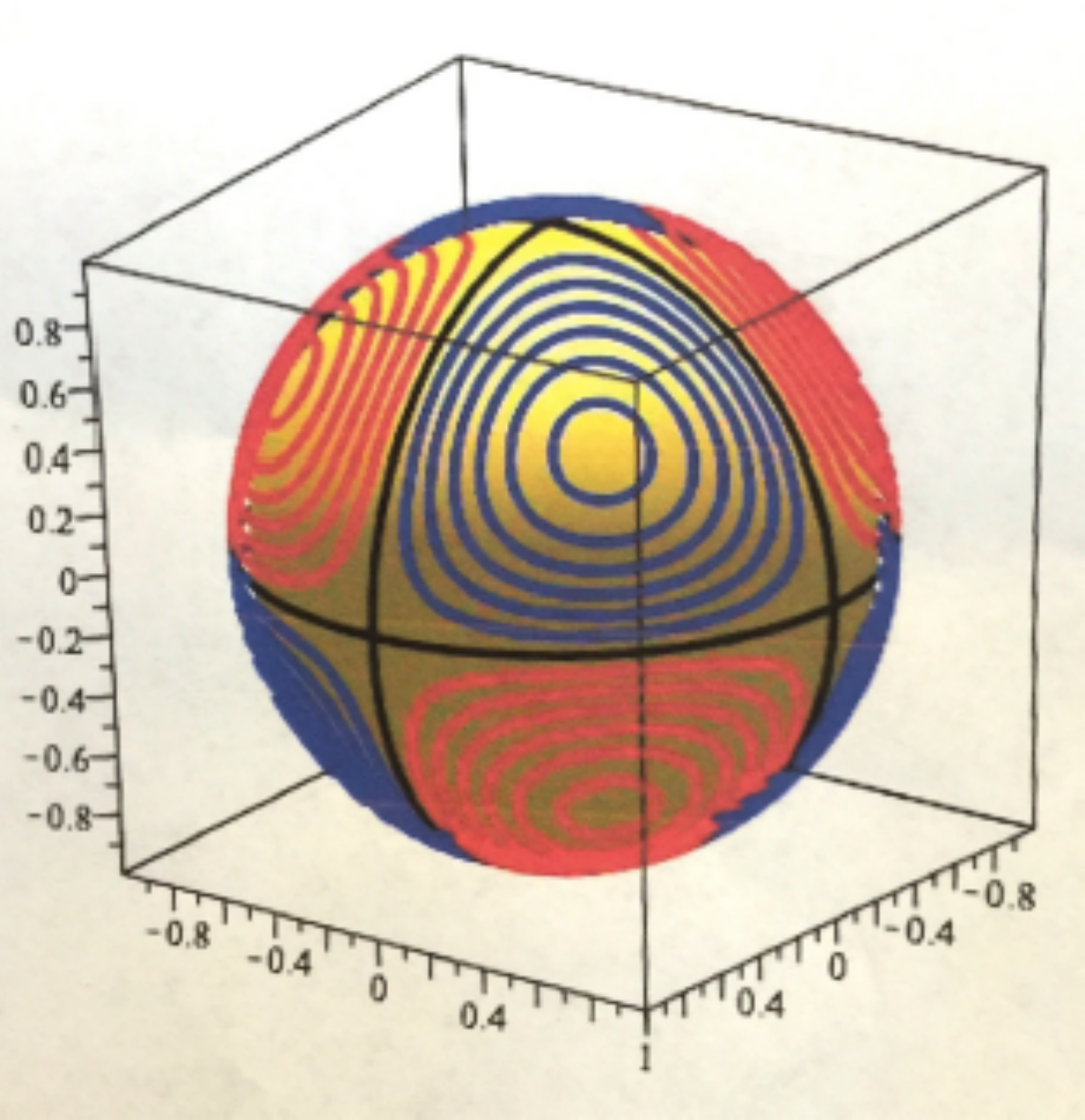}
\caption{{\textit{\textbf{\boldmath Level curves of $\det$ on the $2$-sphere of diagonal matrices}}}}
}
\end{figure}

\vfill
\eject

\noindent\textbf{\boldmath (7) The 7-dimensional variety  $V^7$  of singular matrices on  $S^8(\sqrt{3})$}.
The singular $3\times 3$ matrices  $A$  on  $S^8(\sqrt{3})$  fill out a 7-dimensional algebraic
variety  $V^7$  defined by the equations
$\Vert A \Vert^2  =  3\quad \mbox{and}\quad  \det A  =  0$.
Nothing in our warmup with $2\times 2$ matrices prepares us for the incredible richness 
in the geometry and topology of this variety, which is sketched in Figure 10.

At the lower left is the 4-manifold $M^4$ of matrices of rank 1, along which $V^7$ fails
to be a manifold, and at the upper right is the 5-manifold $M^5$ of best matrices of rank 2 .

The little torus linking $M^4$ signals (in advance of proof) 
that a torus's worth of geodesics on $V^7$ shoot out
orthogonally from each of its points, while the little circle linking $M^5$ signals that a
circle's worth of geodesics on $V^7$ shoot out orthogonally from each of its points.

These are the same geodesics, each an eighth of a great circle, and they fill  $V^7$ with
no overlap along their interiors.

\begin{figure}[h!]
\center{\includegraphics[height=180pt]{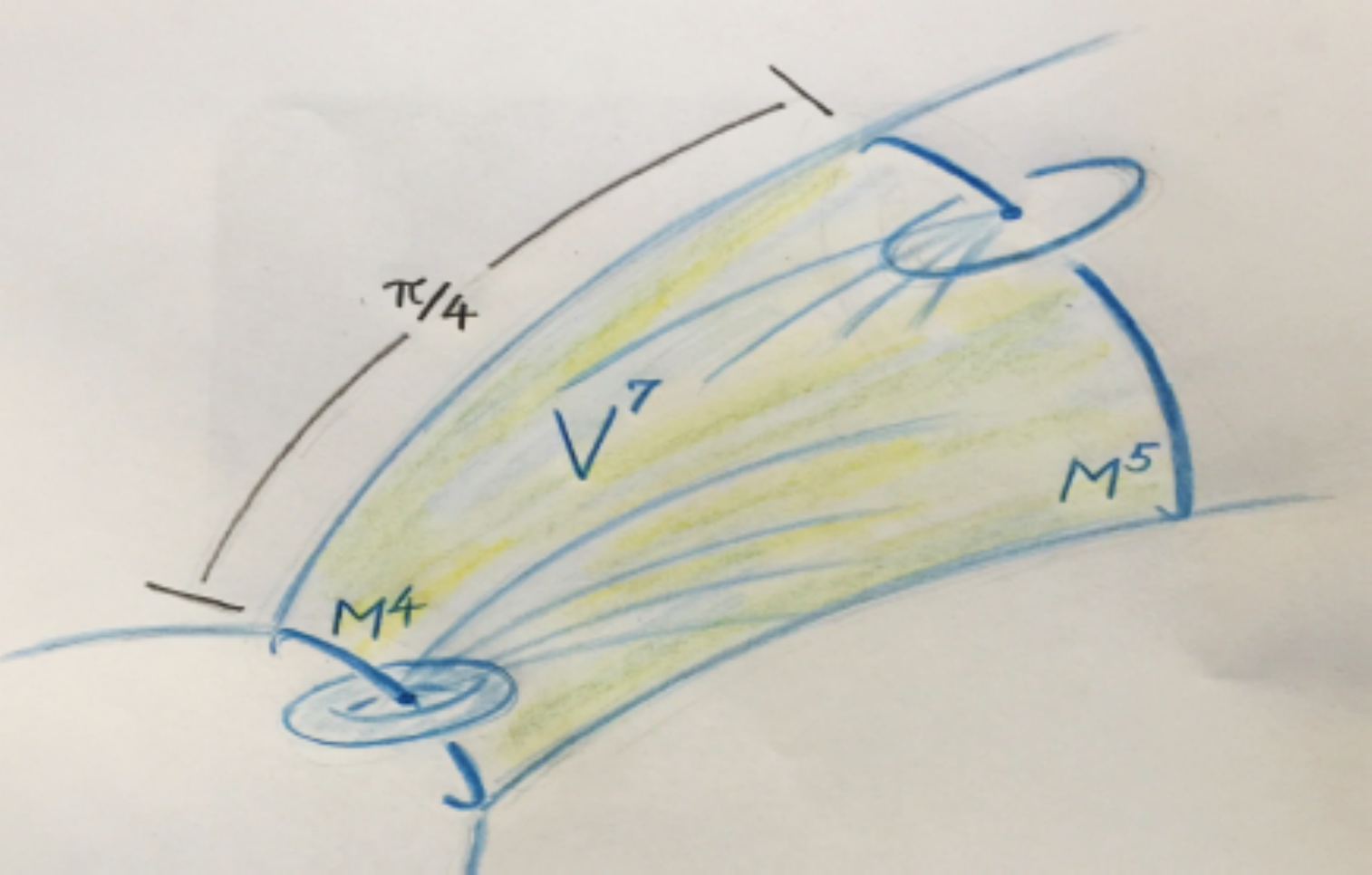}
\caption{{\textit{\textbf{\boldmath The variety $V^7$ of singular $3\times 3$ matrices}}}}
}
\end{figure}

\medskip

\noindent{\textbf{\boldmath (8) What portion of  $V^7$  is a manifold?}
Identifying the set of all $3\times 3$ matrices with Euclidean space  $\reals^9$,  we consider the determinant function
$\det\colon\reals^9\to\reals$.

Let  $A  =  (a_{rs})$  be a given $3\times 3$ matrix.  Then one easily computes the gradient of the determinant function to be 
$$(\nabla\det)_A  =  \sum_{r,s} A_{rs}\frac{\partial}{\partial a_{rs}},$$
where  $A_{rs}$  is the cofactor of  $a_{rs}$  in  $A$.

Thus  $(\nabla\det)_A$  vanishes if and only if all the $2\times 2$ cofactors of  $A$  vanish,
which happens only when  $A$  has rank  $\le  1$.

The subvariety  $V^7$  of  $S^8(\sqrt{3})$  consisting of the singular matrices is 
the zero set of the determinant function $\det\colon S^8(\sqrt{3})\to\reals$.

If $A$  is a matrix of rank 2 on $V^7$ then $\det A=0$ and the gradient vector  $(\nabla\det)_A$  is nonzero there, 
when  det  is considered
as a function from  $\reals^9\to\reals$.
Since $\det (tA)  =  0$  for all real numbers $t$,  the vector  $(\nabla\det)_A$
must be orthogonal to the ray through $A$,  and hence tangent to  $S^8(\sqrt{3})$.
Therefore  $(\nabla\det)_A$  is also nonzero when  det  is considered as a
function from  $S^8(\sqrt{3})\to\reals$.
It follows
that  $V^7$  is a submanifold of  $S^8(\sqrt{3})$  at all its points  $A$  of rank 2.

But $V^7$  fails to be a manifold at all its points
of rank 1,  that is, along the subset  $M^4$, as we will confirm shortly.

\medskip

\noindent\textbf{\boldmath (9) The 5-cell fibres of $N$ and $N'$ meet $V^7$ orthogonally along their 
boundaries.}
This was noted earlier for $2\times 2$ matrices on $S^3(\sqrt{2})$. 

\medskip

\noindent\textbf{Lemma 1. \textit{\boldmath On $S^8(\sqrt{3})$, the gradient vector 
field of the determinant function, when
evaluated at a diagonal matrix, is tangent to the great\break $2$-sphere of diagonal matrices.}}

\noindent\textit{Proof}.\ \  Looking once again at the gradient of the determinant function,
$$(\nabla\det)_A = \sum_{r,s}A_{rs}\frac{\partial}{\partial a_{rs}}\,,$$
where $A_{rs}$ is the cofactor of $a_{rs}$ in $A$, we see that if $A$ is a diagonal matrix,
then $(\nabla\det)_A$ is tangent to the space of diagonal matrices because each off-diagonal
cofactor is zero, and if $A$ lies on 
$S^8(\sqrt{3})$,
then the projection of $(\nabla\det)_A$ to $S^8(\sqrt{3})$ is still tangent to the space of diagonal
matrices there.

\medskip

\noindent\textbf{Lemma 2. \textit{\boldmath More generally, this gradient field is tangent to the 
$5$-dimensional cross-sectional
cells of the tubular neighborhoods $N$ and $N'$ of $SO(3)$ and $O^-(3)$.}}

\medskip

\noindent\textit{Proof}.\ \ If $D^2$ is the 2-dimensional cell of diagonal matrices in the 
tubular neighborhood
$N$ of $SO(3)$ on $S^8(\sqrt{3})$ , then its isometric images $U D^2U^{-1}$, as $U$ ranges over
$SO(3)$, fill out the cross-sectional 5-cell $D^5$ of $N$ at the identity $I$. Since the
determinant function is invariant under this conjugation, its gradient is equivariant,
and so must be tangent to this $D^5$ at each of its points. Then, using left translation
by elements of $SO(3)$, we see that the gradient field is tangent to \textbf{\textit{all}} the cross-sectional
5-cells of $N$ \dots and likewise for $N'$.

\medskip

\noindent\textbf{Proposition 3. \textit{\boldmath The cross-sectional $5$-cell 
fibres of the tubular neighborhoods
$N$ and $N'$ of $SO(3)$ and $O^-(3)$ on $S^8(\sqrt{3})$ meet the variety $V^7$ of singular
matrices orthogonally at their boundaries.}}

\medskip

\noindent\textit{Proof}.\ \ The gradient vector field of the determinant function 
on $S^8(\sqrt{3})$ is orthogonal to the
level surface $V^7$ of this function, and at the same time it is tangent to the cross-sectional
5-cell fibres of $N$ and $N'$. So it follows that these 5-cell fibres meet $V^7$
orthogonally at their boundaries.

\medskip

\noindent\textbf{\boldmath (10) The submanifold  $M^4$  of matrices of rank 1}.
First we identify $M^4$  as a manifold. Define  $S^2\otimes S^2$  to be the quotient of  $S^2 \times S^2$  by the 
equivalence relation
$(x, y ) \sim (-x, -y)$,  a space which is (coincidentally) also homeomorphic to the 
Grassmann manifold of unoriented 2-planes through the origin in real 4-space.
It is straightforward to confirm that $M^4$  is homeomorphic to  $S^2\otimes S^2$.
  
Define a map  $f\colon S^2 \times S^2\to M^4$  by sending the pair of points  $\x  =  (x_1 , x_2 , x_3)$  and  
$\y  =  (y_1 , y_2 , y_3)$  on  $S^2 \times S^2$  to the $3 \times 3$ matrix  $(x_r y_s)$,  scaled up to lie on  $S^8(\sqrt{3})$.
Then check that this map is onto, and that the only duplication is that  $(\x, \y)$  and  $(-\x, -\y)$ go to the same matrix.

\medskip

\noindent\textbf{\textit{Remarks.}}  (1)  $M^4$  is an \textit{\textbf{orientable}} manifold because the involution 
$(\x, \y)\to (-\x, -\y)$  
of   $S^2 \times S^2$  is orientation-preserving.

\noindent (2)  $M^4$  is a single orbit of the $O(3) \times O(3)$  action.

\noindent (3) The integer homology groups of  $M^4$  are
$$H_0(M^4) = \ints , \quad H_1(M^4) = \ints_2 , \quad  H_2(M^4) = \ints_2 ,\quad  H_3(M^4) = 0 \quad\mbox{and}\quad  
H_4(M^4) = \ints,$$
an exercise in using Euler characteristic and Poincar\'e duality (Hatcher [2002]).
Thus $M^4$ has the same \textbf{\textit{rational}} homology as the 4-sphere $S^4$.

\medskip

\noindent\textbf{\boldmath (11) Tangent and normal vectors to  $M^4$}.
At the point  $P  =  \diag (\sqrt{3}, 0, 0)$,  the tangent and normal spaces to $M^4$ within $S^8(\sqrt{3})$ are
$$T_PM^4=\left\{\left[\begin{array}{ccc}0&a&b\\c&0&0\\d&0&0\end{array}\right]\ a,b,c,d\in\reals\right\} \quad
\mbox{and}\quad (T_PM^4)^\perp=
\left\{\left[\begin{array}{ccc}0&0&0\\0&a&b\\0&c&d\end{array}\right]\ a,b,c,d\in\reals\right\}$$
We leave this to the interested reader to confirm.

\medskip

\noindent\textbf{\boldmath(12) The singularity of  $V^7$  along  $M^4$}.
Let  
$A  =\left[\begin{array}{cc}a&b\\c&d\end{array}\right]$
be a $2\times 2$ matrix with  $a^2 + b^2 + c^2 + d^2  =  1$.
Then a geodesic (i.e., great circle) $\gamma(t)$  on  $S^8(1)$  which runs through the 
rank 1 matrix  $P  =  \diag(1, 0, 0)$  at time  $t  =  0$,  and is orthogonal there 
to  $M^4$  has the form
$$\gamma(t)=\left[\begin{array}{ccc}\cos t&0&0\\0&a\sin t&b\sin t\\0&c\sin t&d\sin t\end{array}\right],\quad
\mbox{with} \quad
\gamma'(0)=\left[\begin{array}{ccc}0&0&0\\0&a&b\\0&c&d\end{array}\right].$$
If the $2\times 2$ matrix $A$  above has rank 2,  then $\gamma(t)$ immediately has rank 3  for small  $t > 0$.  But if  
$A$  has rank 1,  then $\gamma(t)$  has only rank 2  for small  $t > 0$.

\begin{wrapfigure}[13]{r}{0.50\textwidth}
\vspace{-20pt}
\center{\includegraphics[width=0.48\textwidth]{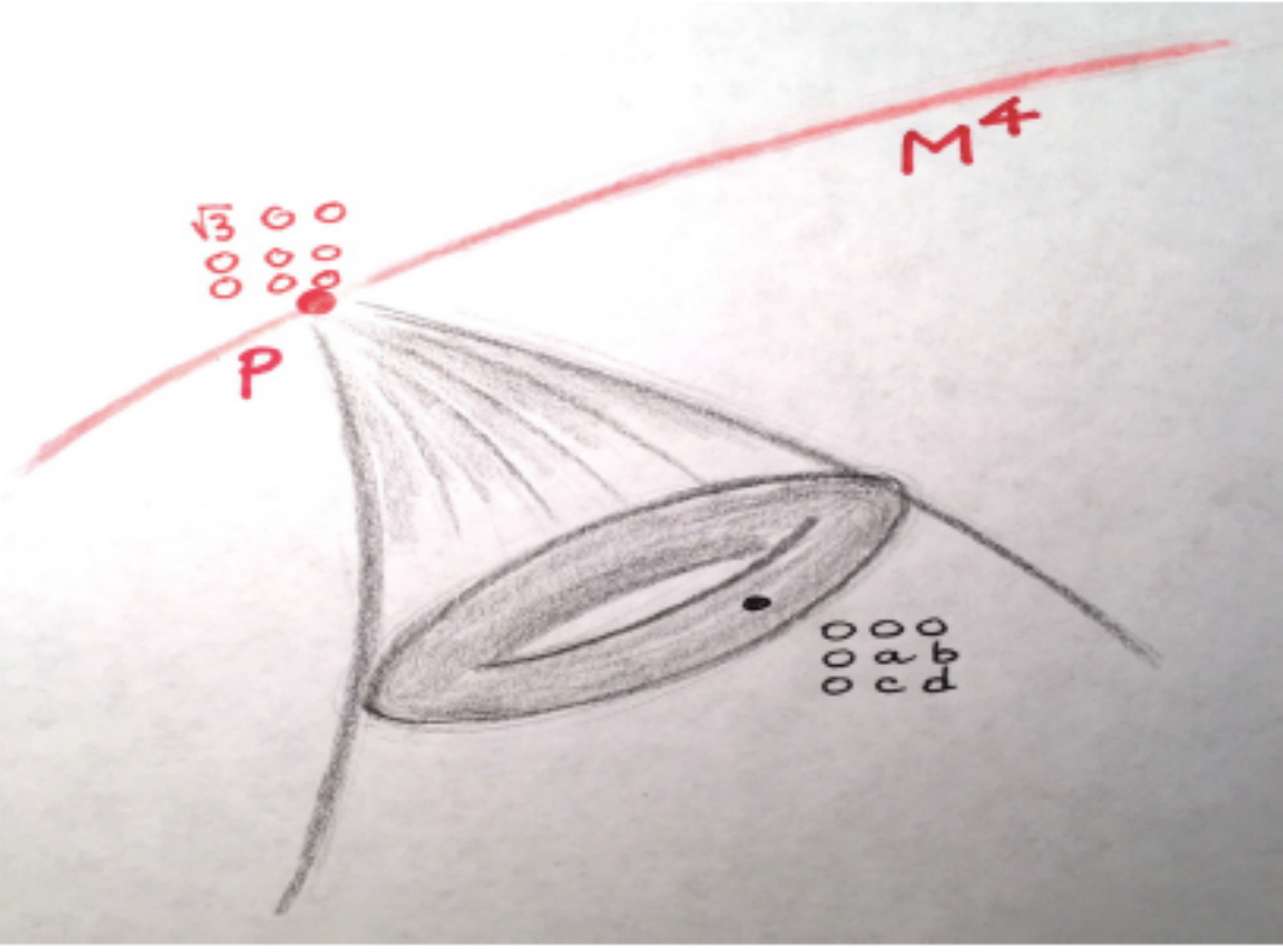}
\caption{{\textit{\textbf{\boldmath The normal cone to $M^4$ in $V^7$ at the point $P$ is a cone over a Clifford torus 
$a^2+b^2+c^2+d^2=1$ and $ad-bc=0$}}}}
}
\end{wrapfigure}
We know from our study of $2\times 2$ matrices that those of rank 1 form a cone (punctured at the origin) over the 
Clifford torus in  $S^3(1)$. Thus the tubular neighborhood of  $M^4$  in  $V^7$  is a bundle over  $M^4$  whose normal fibre is a cone over the Clifford torus.
We indicate this pictorially in Figure 11.

One can use the information above to show that  

\noindent (1) $V^7 - M^5$  is a ``tubular''
neighborhood of  $M^4$,  whose cross-sections are great circle cones of
angular radius $\pi/4$ over Clifford tori on great 3-spheres which meet  $M^4$
orthogonally.

\noindent (2) The 2-torus's worth of geodesic rays shooting out from each point 
of  $M^4$  in  $V^7$  terminate along a full 2-torus's worth of points in  $M^5$. 

\medskip

\noindent\textbf{\boldmath(13) The submanifold  $M^5  =  \{\mbox{Best of rank 2}\}$.}
Recall that the ``best'' $3 \times 3$ matrices of rank 2 are 
those which are orthogonal on a
2-plane through the origin, and zero on its orthogonal complement.

An example of such a matrix is  
$P  =  \diag(1, 1, 0) =\left[\begin{array}{ccc}1&0&0\\0&1&0\\0&0&0\end{array}\right]$, 
representing orthogonal projection of  $x y z$-space to the  $x y$-plane.

We let  $M^5$  denote the set of best $3 \times 3$ matrices of rank 2,  scaled up to lie 
on  $S^8(\sqrt{3})$.  This set  is a single orbit of the  $SO(3)\times SO(3)$  action on $\reals^9$.

\medskip

\noindent\textbf{Claim: \textit{\boldmath  $M^5$  is homeomorphic to  $\reals P^2 \times  \reals P^3$.}}

\medskip

\noindent\textit{Proof.}\ \   Let  $T$  be one of these best $3\times 3$ matrices of rank 2 .
Then the kernel of  $T$  is some unoriented line through the origin in  $\reals^3$,
hence an element of  $\reals P^2$.

An orthogonal transformation of  $(\ker T)^\perp$  to a 2-plane through the origin 
in  $\reals^3$  can be uniquely extended to an orientation-preserving orthogonal transformation  
$A_T$  of $\reals^3$  to itself, hence an element of $SO(3)$.

Then the correspondence  $T \to (\ker T\,,\, A_T)$  gives the homeomorphism of 
$M^5$  with  $\reals P^2  \times  SO(3)$,  equivalently, with  $\reals P^2  \times  \reals P^3$.

\medskip

\noindent\textit{\textbf{Remark.}}  $M^5$ is non-orientable, and its integer homology groups  are
$$H_0(M^5)  =  \ints ,\quad  H_1(M^5)  =  \ints_2 + \ints_2 ,\quad  H_2(M^5)  =  \ints_2 , \quad H_3(M^5)  =  \ints + \ints_2 ,$$  
$$H_4(M^5)  =  \ints_2 ,\quad  H_5(M^5)  =  0,$$
an exercise in using the K\"unneth formula (Hatcher [2002]) for the homology of a product.

\medskip

\noindent\textbf{\boldmath(14) Tangent and normal vectors to  $M^5$}. 
At the point  $P  =  \diag \left(\sqrt{\frac32}\,,\,\sqrt{\frac32}\,,\, 0\right)$,  the tangent and normal spaces to $M^5$ within 
$V^7$ are
$$T_PM^5=\left\{\left[\begin{array}{ccc}0&-a&b\\a&0&c\\d&e&0\end{array}\right]\ a,b,c,d,e\in\reals\right\} \quad
\mbox{and}\quad (T_PM^5)^\perp=
\left\{\left[\begin{array}{ccc}a&b&0\\b&-a&0\\0&0&0\end{array}\right]\ a,b\in\reals\right\}$$
and $(T_PV^7)^\perp\subset T_PS^8(\sqrt{3})$ is spanned by $\diag(0,0,1)$, as the reader can confirm.

\medskip

\noindent\textbf{\boldmath(15) The tubular neighborhood of  $M^5$  inside  $V^7$ }.

\medskip

\noindent\textbf{Claim: \textit{\boldmath  $V^7 -  M^4$  is a tubular neighborhood of  $M^5$,  whose cross sections are 
round cells of angular radius  $\pi/4$  on great $2$-spheres which meet  $M^5$  orthogonally.}}

\medskip

\noindent\textit{Proof.}\ \  We start on  $M^5$  at the scaled point  
$P  =  \diag(1, 1, 0) =\left[\begin{array}{ccc}1&0&0\\0&1&0\\0&0&0\end{array}\right]$,
which represents orthogonal projection of  $x y z$-space to the  $x y$-plane.
Then we consider the tangent vectors
$$T_1=\left[\begin{array}{ccc}0&1&0\\1&0&0\\0&0&0\end{array}\right]
\quad\mbox{and}\quad
T_2=\left[\begin{array}{ccc}1&0&0\\0&-1&0\\0&0&0\end{array}\right]$$
which are an orthogonal basis for $(T_PM^5)^\perp\subset  T_PV^7$.

If we exponentiate the vector in  $(T_PM^5)^\perp$ given by  $a T_1  +  b T_2$,  with  $a^2 + b^2  =  1$,  
from the point  $P$,  we get
$$P(t)  =  (\cos t) P  +  (\sin t) (a T_1  +  b T_2)  =
\left[\begin{array}{ccc}\cos t+b\sin t&a\sin t&0\\a\sin t&\cos t-b\sin t&0\\0&0&0\end{array}\right]$$
which has rank 2 for  $0 \le t  <\pi/4$.   All these matrices have the same 
kernel and same image as  $P$ .
But   
$\displaystyle P(\pi/4)  =  \frac{1}{\sqrt{2}}\left[\begin{array}{ccc}1+b&a&0\\a&1-b&0\\0&0&0\end{array}\right]$,
which only has rank 1, and therefore lies on  $M^4$.

The set of points  $\{P(t): 0 \le  t \le\pi/4\}$  is one-eighth of a great circle on  $S^8(\sqrt{3})$,  beginning at the point  $P  =  P(0)$  
on  $M^5$  and ending at the point  $P(\pi/4)$  on  $M^4$.  
Let's call this set a \textit{\textbf{ray}}.  

We see from the entries in the above matrix that the circle's worth of rays shooting 
out from the point  $P$  orthogonal to  $M^5$  in  $V^7$  terminate along a full circle's worth 
of points on  $M^4$.  We can think of this as an ``absence of focusing''.

Since  $M^5$  is a single orbit of the  $SO(3) \times SO(3)$  action on  $\reals^9$,  the above situation 
at the point  $P$  on  $M^5$  is replicated at every point of  $M^5$,  confirming the claim made above.  

\medskip

\noindent\textbf{\boldmath(16) The wedge norm on  $V^7$.} Recall that for 
$2 \times 2$ matrices viewed as points in $\reals^4$  and then restricted to  $S^3(\sqrt{2})$,
the determinant function varies between a maximum of  1  on  $SO(2)$  and a minimum of $-1$  on  $O^-(2)$,  with the middle 
value zero assumed on the Clifford torus of singular matrices.  The level sets of this for values strictly between 
$-1$  and  1  are tori parallel to the Clifford torus, and are 
principal orbits of the  $SO(2) \times SO(2)$ action. Their orthogonal trajectories are the geodesic arcs 
leaving  $SO(2)$  orthogonally and arriving at  $O^-(2)$  orthogonally a quarter of a great circle later.

We seek a corresponding function on the variety  $V^7$  of singular matrices on $S^8(\sqrt{3})$, 
whose level sets fill the space between $M^4$ and $M^5$, 
and
to this end, turn to the wedge norm  $\Vert A\wedge A\Vert$,  defined as follows.

If  $A\colon V\to W$  is a linear map between the real vector spaces  $V$  and  $W$,  then the induced linear map  
$A\wedge A \colon \wedge^2V\to\wedge^2 W$  between spaces of 2-vectors is defined by
$$(A\wedge A) (\v_1\wedge \v_2)  =  A(\v_1)\wedge A(\v_2),$$
with extension by linearity. If  $V  =  W  = \reals^2$,  then the space $\wedge^2\reals^2$  is one-dimensional,
and  $A\wedge A$  is simply multiplication by  $\det A$, while if
$V  =  W  =  \reals^3$,  then the space $\wedge^2\reals^3$  is three-dimensional,
and  $A\wedge A$  coincides with the matrix  of cofactors of  $A$.

The wedge norm is defined by
$\Vert A\wedge A\Vert^2=\sum_{i,j}(A\wedge A)_{ij}^2$,
and 
is  easily seen to be $SO(3) \times SO(3)$-invariant, and thus
constant along the orbits of this action.
It has the following properties:
\begin{enumerate}
\item[(1)] On $V^7$ the wedge norm takes its maximum value of $3/2$ on $M^5$ and
its minimum value of 0 on $M^4$.
\item[(2)] The level sets between these two extreme values are 6-dimensional
submanifolds which are principal orbits
of the $SO(3) \times SO(3)$ action.
\item[(3)] The orthogonal trajectories of these level sets are geodesic arcs,
each an eighth of a great circle, meeting both $M^4$ and $M^5$ orthogonally.
\end{enumerate}

\medskip

\noindent\textbf{\boldmath (17) Concrete generators for the 4-dimensional homology of $V^7$.}
If we remove both components of the orthogonal group $O(3)$ from the 8-sphere
$S^8(\sqrt{3})$, then what is left over deformation retracts to the variety $V^7$, since each
cross-sectional 5-cell $D^5$ in the tubular neighborhoods of these two components
has now had its center removed, and so can deformation retract to its boundary
along great circle arcs.

Therefore $V^7$ has the same integer homology as $S^8(\sqrt{3}) - O(3)$, which can be
computed by Alexander duality (Hatcher [2002]), and we learn that
$$H_0(V^7) = \ints, \quad H_4(V^7) = \ints + \ints,\quad H_5(V^7) = \ints_2 + \ints_2 ,\quad H_7(V^7) = \ints,$$
while the remaining homology groups are zero.
The variety $V^7$ is orientable because it divides $S^8(\sqrt{3})$ into two components,
but its homology is excused from satisfying Poincar\'e duality because it is
not a manifold.

We seek concrete cycles generating $H_4(V^7)$.

Pick a point on each component of $O(3)$, for example, the identity $I$ on $SO(3)$,
and $-I$ on $O^-(3)$.
Then go out a short distance in the cross-sectional 5-cells of the two tubular
neighborhoods, and we will have a pair of 4-spheres, each linking the corresponding
component of $O(3)$, and therefore generating $H_4(S^8(\sqrt{3}) - O(3))$.
Pushing these 4-spheres outwards to $V^7$ along the great circle rays of these two
5-cells provides the desired generators for $H_4(V^7)$.

How are these generators positioned on $V^7$?

\begin{wrapfigure}[15]{r}{0.5\textwidth}
\vspace{-35pt}
\center{\includegraphics[width=0.48\textwidth]{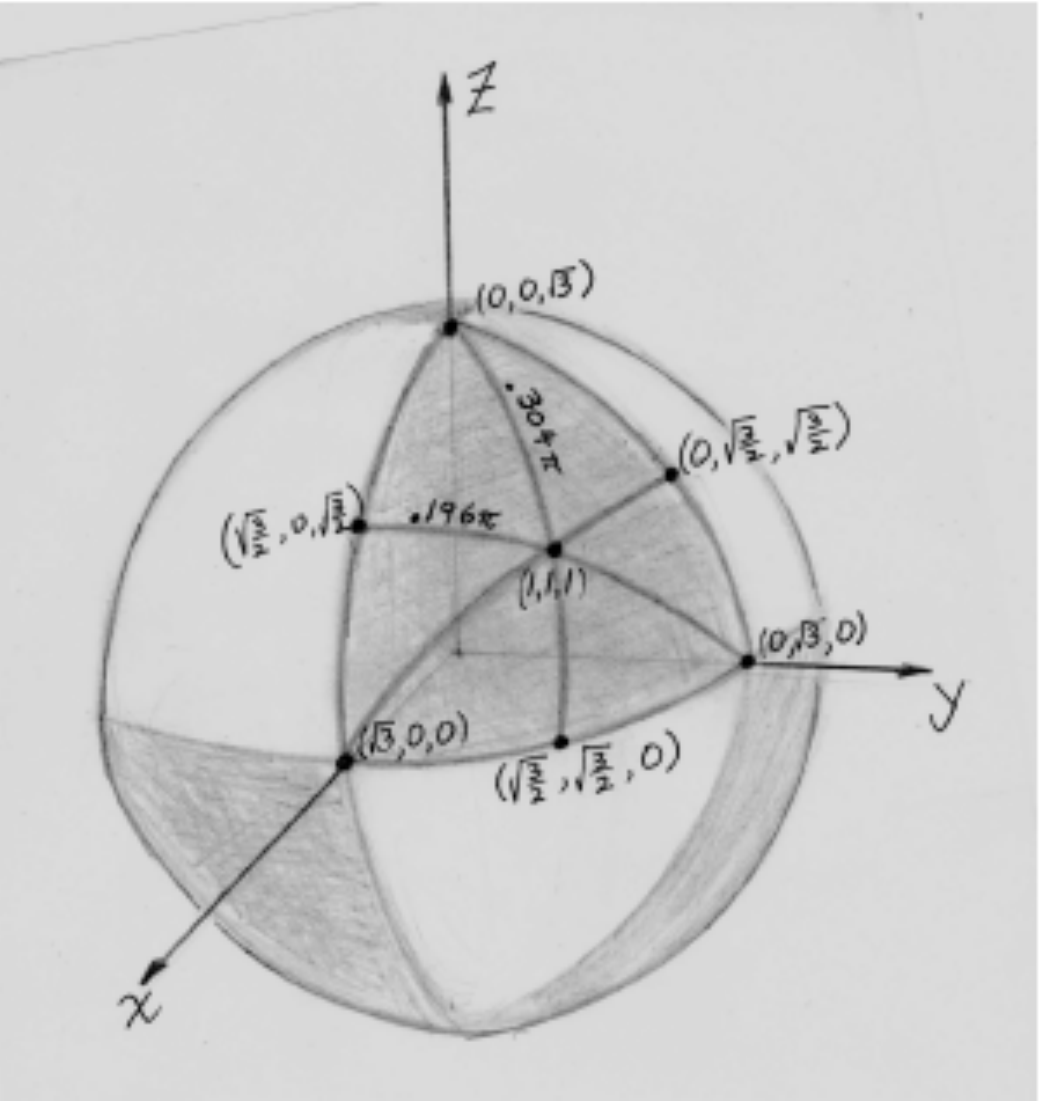}
\caption{{\textit{\textbf{\boldmath Diagonal matrices in $S^8(\sqrt{3})$}}}}
}
\end{wrapfigure}
The key to the answer can be found in the diagonal $3 \times 3$ matrices.
In Figure 12, consider the spherical triangle centered at $(1, 1, 1)$.
We noted earlier that the three vertices of this triangle lie in $M^4$, and the centers
of its three edges in $M^5$. The six half-edges are geodesics, each an eighth of a
great circle.

Conjugating by $SO(3)$ promotes this triangle to the cross-sectional 5-cell centered
at the identity in the tubular neighborhood of $SO(3)$, and promotes the decomposition
of the boundary of the triangle to a decomposition of the boundary $+S^4$ of
this 5-cell.

\begin{figure}[h!]
\center{\includegraphics[height=220pt]{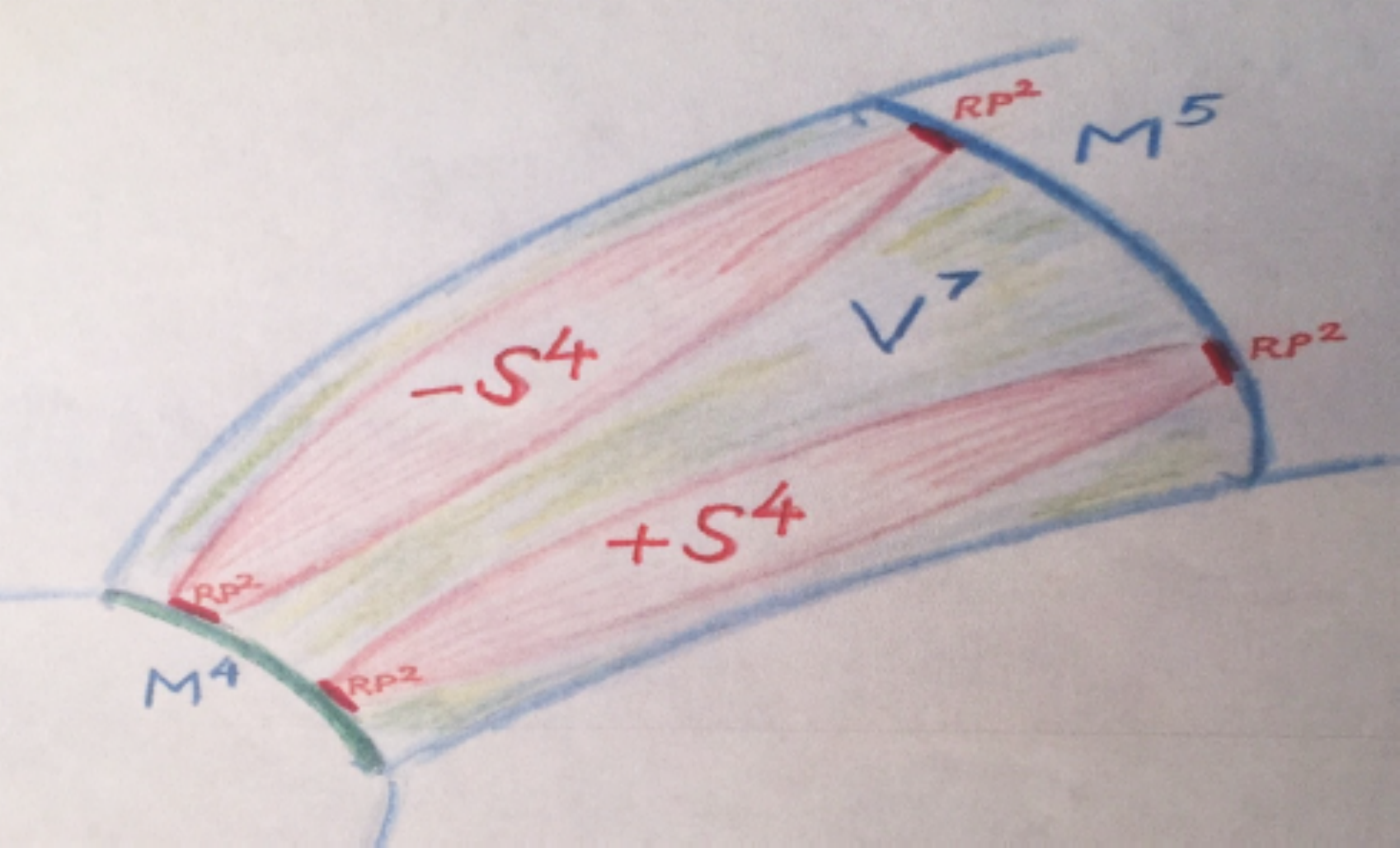}
\caption{{\textit{\textbf{\boldmath Generators of $H_4(V^7)$ are $4$-spheres with $\reals P^2$ ends in $M^4$ and $M^5$}}}}
}
\end{figure}

This is enough to reveal the positions of our two generators of $H_4(V^7)$. We show this in Figure 13, where
\begin{enumerate}
\item[(1)] The lower 4-sphere $+S^4$ links $SO(3)$ in $S^8(\sqrt{3})$, and is the set of
symmetric positive semi-definite matrices there which are not positive definite.
\item[(2)] The upper 4-sphere $-S^4$ links $O^-(3)$, and is the set of symmetric negative semi-definite
matrices on $S^8(\sqrt{3})$ which are not negative definite.
\item[(3)] Each of these 4-spheres has an $\reals P^2$ end in $M^4$ and another $\reals P^2$ end in $M^5$,
and is smooth, except at the end in $M^4$.
\item[(4)] The $SO(3)$ action by conjugation on each 4-sphere is the same as that on the unit
4-sphere in the space of traceless, symmetric $3\times 3$ matrices described by Blaine Lawson
[1980] . The principal orbits are all copies of the group $S^3$ of unit quaternions, modulo
its subgroup $\{\pm1 , \pm i , \pm j , \pm k\}$, the singular orbits are the $\reals P^2$ ends, and the orthogonal
trajectories are geodesic arcs, each an eighth of a great circle.
\end{enumerate}

\vspace{-6pt}

\noindent\textbf{(18)  Nearest orthogonal neighbor.}  Start with a nonsingular $3\times 3$ matrix  $A$  on  
$S^8(\sqrt{3})$  and suppose, to be specific, that  $A$  lies in the open tubular neighborhood $N$ 
of  $SO(3)$. 

 \textbf{\textit{\boldmath We claim that the nearest orthogonal neighbor to  $A$  on that $8$-sphere 
is the center of the cross-sectional fibre of  $N$  on which it lies.}}

To see this, note that a geodesic (great circle arc) from  $A$  to its nearest neighbor  $U$
on  $SO(3)$ must meet  $SO(3)$  orthogonally at  $U$,  and therefore must lie in the cross-sectional fibre of  $N$  
through  $U$.  It follows that  $A$  also lies in that fibre, whose center 
is at  $U$,  confirming the above claim.

\medskip

\vspace{-6pt}

\noindent\textbf{(19)  Nearest singular neighbor.} Start with a nonsingular $3\times 3$ matrix $A$ on 
$S^8(\sqrt{3})$, say with $\det A > 0$.

\textit{\textbf{\boldmath We claim that the nearest singular neighbor to $A$ on that 
$8$-sphere lies on the
boundary of the cross-sectional $5$-cell of the tubular neighborhood $N$ of $SO(3)$
which contains $A$.}}

\medskip

\noindent\textbf{Lemma. \textit{\boldmath Let $A$ be a nonsingular $3\times 3$ matrix on 
$S^8(\sqrt{3})$, and let $B$ be the
closest singular matrix to $A$ on this 8-sphere. Then $B$ has rank $2$.}}

\medskip

\noindent\textit{Proof.}\ \  Suppose $B$ has rank 1, and therefore lies in $M^4$. 
Since $SO(3) \times SO(3)$
acts transitively on $M^4$, we can choose orthogonal matrices $U$ and $V$ so that
$U B V^{-1} = \diag (\sqrt{3}, 0, 0)$, and this will then be the closest singular matrix to the
nonsingular matrix $U A V^{-1}$.
So we can assume that $B = \diag(\sqrt{3}, 0, 0)$ already.

Since $A$ is nonsingular, it must have at least one nonzero entry $a_{ij}$ for some
$i > 1$ and $j > 1$. Now let $T$ be the matrix with all zeros except in the $ij$th spot,
with $t_{ij} = \sgn(a_{ij}) \sqrt{3}$. Then $T$ also lies on $S^8(\sqrt{3})$ and is orthogonal to $B$.

Hence the matrices $B(t) = \cos t\, B + \sin t\, T$ lie on $S^8(\sqrt{3})$ as well, and
$$\langle A , B(t) \rangle = \cos t\, \langle A , B \rangle + \sin t\, \langle A , T \rangle.$$
The derivative of this inner product with respect to $t$ at $t = 0$ is
$$\langle A , T \rangle = \vert a_{ij}\vert\sqrt{3}>0.$$
Therefore, for small values of $t$, $B(t)$ is a matrix of rank 2 on $S^8(\sqrt{3})$ that is
closer to $A$ than $B$ was. This contradicts the assumption that $B$ was closest
to $A$, and proves the lemma.

\noindent\textbf{\textit{Remarks.}} (1) For visual evidence in support of this lemma, 
look at the front shaded
spherical triangle on the great 2-sphere of diagonal $3 \times 3$ matrices in 
Figure 12, and
note that if $A$ is an interior point of this triangle, then the closest boundary point $B$
cannot be one of the vertices.

\noindent(2) More generally, let $A$ be an $n \times n$ matrix of rank $> r$ on $S^{n^2 -1}$. 
Then the matrix $B$ on $S^{n^2 -1}$ of rank $\le r$ that is closest to $A$ actually has rank $r$.

\medskip

Now given the nonsingular matrix $A$ on $S^8(\sqrt{3})$, its nearest neighbor $B$ on $V^7$
must have rank 2, and therefore lies in the manifold portion of $V^7$. It follows
that the shortest geodesic from $A$ to $B$ is orthogonal to $V^7$, and since we saw
in section 9 that the 5-cell fibres of $N$ meet $V^7$ orthogonally along their boundaries,
this geodesic must lie in the 5-cell fibre of $N$ containing $A$.

Therefore $B$ lies on the boundary of this 5-cell fibre, as claimed.

\noindent\textit{\textbf{Remark.}} Because the 5-cell fibres of $N$ 
are not round, the nearest singular neighbor $B$
is typically \textit{not} at the end of the ray from the center of the cell through $A$, as was true
for $2\times 2$ matrices. We will shortly state the classical theorem of Eckart and Young
which describes this nearest singular neighbor explicitly in terms of singular values.

\bigskip

\centerline{\LARGE{\textbf{\textit{Matrix Decompositions}}}}

\bigskip

\medskip

\centerline{\Large{\textbf{\boldmath{Singular value decomposition}}}}

\medskip

\vspace{-10pt}

Let $A$ be an $n \times k$ matrix, thus representing a linear map $A\colon \reals^k\to\reals^n$.

We seek a matrix decomposition of $A$,
$$A = W D V^{-1},$$
where $V$ is a $k \times k$ orthogonal matrix, where $D$ is an $n \times k$ diagonal matrix,
$$D = \diag(d_1 , d_2 ,\ldots, d_r) , \quad\mbox{with}\  d_1 \ge d_2 \ge \cdots  \ge d_r \ge 0,$$
with $r = \min(k, n)$, and where $W$ is an $n \times n$ orthogonal matrix.

\begin{wrapfigure}[11]{r}{0.5\textwidth}
\vspace{-20pt}
\center{\includegraphics[width=0.48\textwidth]{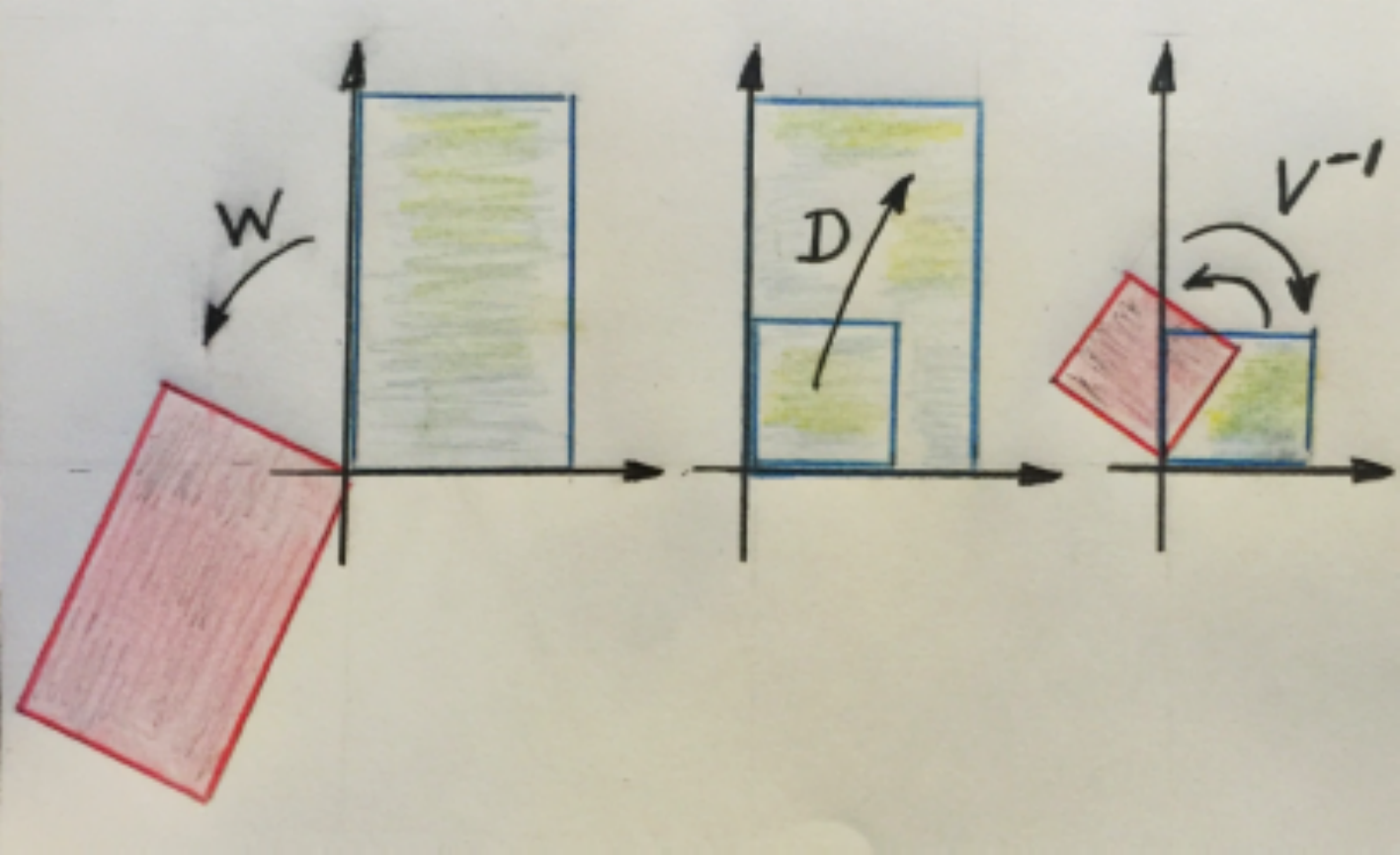}
\caption{{\textit{\textbf{\boldmath Singular value decomposition: $A=WDV^{-1}$}}}}
}
\end{wrapfigure}
The message of this decomposition is that $A$ takes some right angled\break $k$-dimensional
box in $\reals^k$ to some right angled box of dimension $\le k$ in $\reals^n$, with the columns of the
orthogonal
matrices $V$ and $W$ serving to locate the edges of the domain and image boxes, and the diagonal
matrix $D$ reporting expansion and compression of these edges (Figure 14).
See Golub and Van Loan [1996] and Horn and Johnson [1991] for derivation of
this singular value decomposition.

\noindent\textbf{Remarks}.\ \ (1) 
Consider the map $A^T A\colon\reals^k\to\reals^k$, and note that
$$A^T A = (V D W^{-1}) (W D V^{-1}) = V D^2 V^{-1},$$
with eigenvalues $d_1^2,d_2^2,\ldots,d_r^2$ and if $r=n<k$, then also with $k-n$ zero
eigenvalues. 
The orthonormal columns $\v_1, \v_2 ,\ldots, \v_k$ of $V$ are
the corresponding eigenvectors of $A^T A$, since for example
$$A^T A (\v_1) = V D^2 V^{-1} (\v_1) = V D^2 (1, 0,\ldots, 0) = V (d_1^2,0,\ldots,0)=d_1^2\v_1,$$
and likewise for $\v_2 ,\ldots, \v_k$.

\noindent (2) In similar fashion, consider the map $A A^T\colon\reals^n\to\reals^n$, note that
$$A A^T = (W D V^{-1}) (V D W^{-1}) = W D^2 W^{-1},$$
with eigenvalues $d_1^2,d_2^2,\ldots,d_r^2$, and if $r=k<n$, then also with $n-k$ zero
eigenvalues. 
The orthonormal columns  $\w_1,\w_2 ,\ldots, \w_n$ of $W$ are
the corresponding eigenvectors of $A A^T$.

\noindent (3) The singular value decomposition was discovered independently by the
Italian differential geometer Eugenio Beltrami [1873] and the French algebraist
Camille Jordan [1874a, b] , in response to a question about the bi-orthogonal equivalence of quadratic forms.
Later, Erhard Schmidt [1907] introduced the infinite-dimensional analogue of the singular
value decomposition and
addressed the problem of finding the best approximation
of lower rank to a given bilinear form.

Carl Eckart and Gale Young [1936] extended the singular value decomposition to
rectangular matrices, and rediscovered Schmidt's 1907 theorem about approximating
a matrix by one of lower rank.

\noindent (4) Since finding the singular value decomposition of a matrix $A$ is equivalent 
to computing the eigenvalues and
orthonormal eigenvectors of the symmetric matrices $A^T A$ and $A A^T$, all of the computational techniques
that apply to positive (semi)definite symmetric matrices apply, in particular the
celebrated QR-algorithm, which was proposed
independently by John Francis [1961] and Vera Kublanovskaya [1962].
Its later refinement, the \textit{implicitly shifted QR algorithm}, was named one of the top ten
algorithms of the 20th century by the editors of SIAM news (Cipra [2000]).
For more historical details, we recommend Stewart [1993].

\bigskip

\medskip

\centerline{\Large{\textbf{\boldmath{Polar decomposition}}}}

\medskip

\vspace{-10pt}

The \textbf{\textit{polar decomposition}} of an $n\times n$ matrix $A$ is the factoring
$$A = U P,$$
\begin{wrapfigure}[11]{r}{0.48\textwidth}
\vspace{-30pt}
\center{\includegraphics[width=0.43\textwidth]{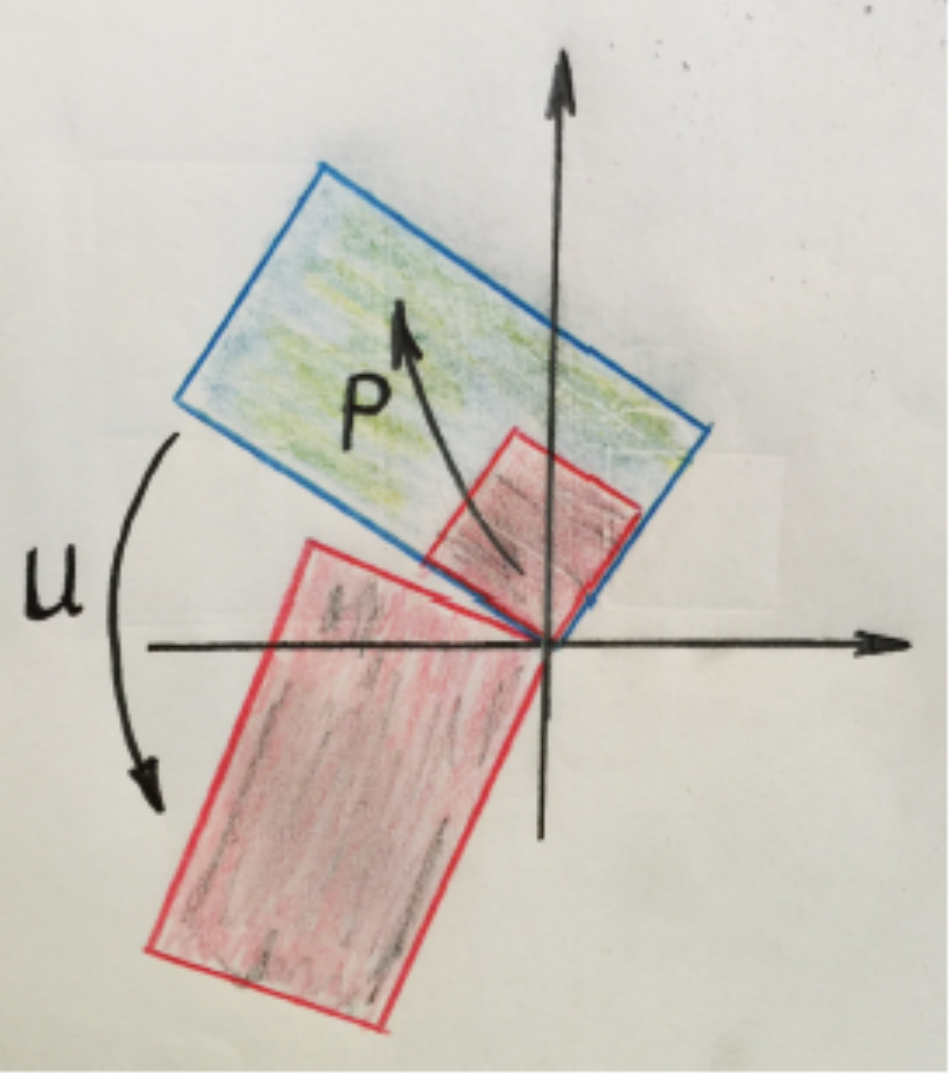}
\caption{{\textit{\textbf{\boldmath Polar decomposition: $A=UP$}}}}
}
\end{wrapfigure}
where $U$ is orthogonal and $P$ is symmetric positive semi-definite.

The message of this decomposition is that $P$ takes some right angled $n$-dimensional
box in $\reals^n$ to itself, edge by edge, expanding and compressing some while perhaps
sending others to zero, after which $U$ moves the image box rigidly to another position (Figure 15).

See Golub and Van Loan [1996] and Horn and Johnson [1991] for derivation of this
polar decomposition.
\vfill
\eject

\noindent\textbf{Remarks}.\ \ (1) Existence of the polar decomposition 
follows immediately from the singular value decomposition for $A$:
$$A = W D V^{-1} = (W V^{-1}) (V D V^{-1}) = U P.$$
Furthermore, if $A = U P$ , then $A^T = P^T U^T = P U^{-1}$, and hence
$$A^T A = (P U^{-1}) (U P) = P^2.$$
Now the symmetric matrix $A^T A$ is positive semi-definite, and 
has a unique symmetric positive semi-definite square root $P = \sqrt{A^TA}$.

\noindent (2) In the polar decomposition $A=UP$, the factor $P$
is uniquely determined by $A$, while the factor $U$ is uniquely determined
by $A$ if $A$ is nonsingular, but not in general if $A$ is singular.

\noindent (3) If $n=3$ and $A$ is nonsingular, with polar decomposition $A = U P$, and if we scale $A$
to lie on $S^8(\sqrt{3})$, then $P$ will also lie on that sphere, and the polar decomposition of $A$ is
just the product coordinatization of the open tubular neighborhoods $N$ and $N'$ of
$SO(3)$ and $O^-(3)$.

\noindent (4) An $n\times n$ matrix $A$ of rank $r$ has a factorization $A = U P$,
with $U$ best of rank $r$ and $P$ symmetric positive semi-definite, and with both
factors $U$ and $P$ uniquely determined by $A$ and having the same rank $r$ as $A$.

\noindent (5) Let $A$ be a real nonsingular $n\times n$ matrix, and let $A = U P$
be its polar decomposition. Then $U$ is the nearest orthogonal matrix to $A$, in
the sense of minimizing the norm $\Vert A-V\Vert$ over all orthogonal matrices $V$.

\noindent (6)  Let $A=UP$ be an $n\times n$ matrix of rank $r$ 
with $U$ best of rank $r$ and $P$ symmetric positive semi-definite.
Then $U$ is the nearest best of rank $r$ matrix to $A$,
in the sense of minimizing the norm
 $\Vert A - V \Vert$ over all
best of rank $r$ matrices $V$.

\noindent (7) The decomposition $A = U P$ is called \textit{\textbf{right polar decomposition}}, to
distinguish it from the \textit{\textbf{left polar decomposition}}
$A = P' U' .$
Given the right polar decomposition $A = U P$, we can write
$A = U P = (U P U^{-1}) U = P' U$
to get the left polar decomposition. If $A$ is nonsingular, then the unique orthogonal
factor $U$ is the same for both right and left polar decompositions, but the symmetric
positive semi-definite factors $P$ and $P'$ are not.
No surprise about the orthogonal factor being the same, since in either case it is the
unique element of the orthogonal group $O(n)$ closest to $A$.

\noindent (8) L\'eon Autonne [1902], in his study of matrix groups, first introduced the polar
decomposition $A = U P$ of a square matrix $A$, where $U$ is unitary and
$P$ is Hermitian, and quickly proved its existence.

\bigskip

\centerline{\Large{\textbf{\boldmath{The Nearest Singular Matrix}}}}

\medskip

\noindent\textbf{Theorem (Eckart and Young, 1936). \textit{\boldmath Let $A$ be an $n \times k$ matrix of rank $r$,
with singular value decomposition $A = W D V^{-1}$, where $V$ is a $k \times k$ orthogonal
matrix, where $D$ is an $n \times k$ diagonal matrix,
$$D = \diag (d_1 , d_2 , \ldots, d_r , 0, \ldots, 0) , \quad\mbox{with}\  d_1 \ge d_2 \ge \ldots ³ d_r > 0 ,$$
and where $W$ is an $n \times n$ orthogonal matrix.\\[0.1cm]
Then the nearest $n \times k$ matrix $A'$ of rank $\le r' < r$ is given by $A' = W D' V^{-1}$,
with $W$ and $V$ as above, and with
$$D' = \diag(d_1 , d_2 , \ldots, d_{r'} , 0 ,\ldots, 0).$$
}}
We illustrate this in Figure 16 in the setting of 3 x 3 matrices. 
\begin{figure}[h!]
\center{\includegraphics[height=180pt]{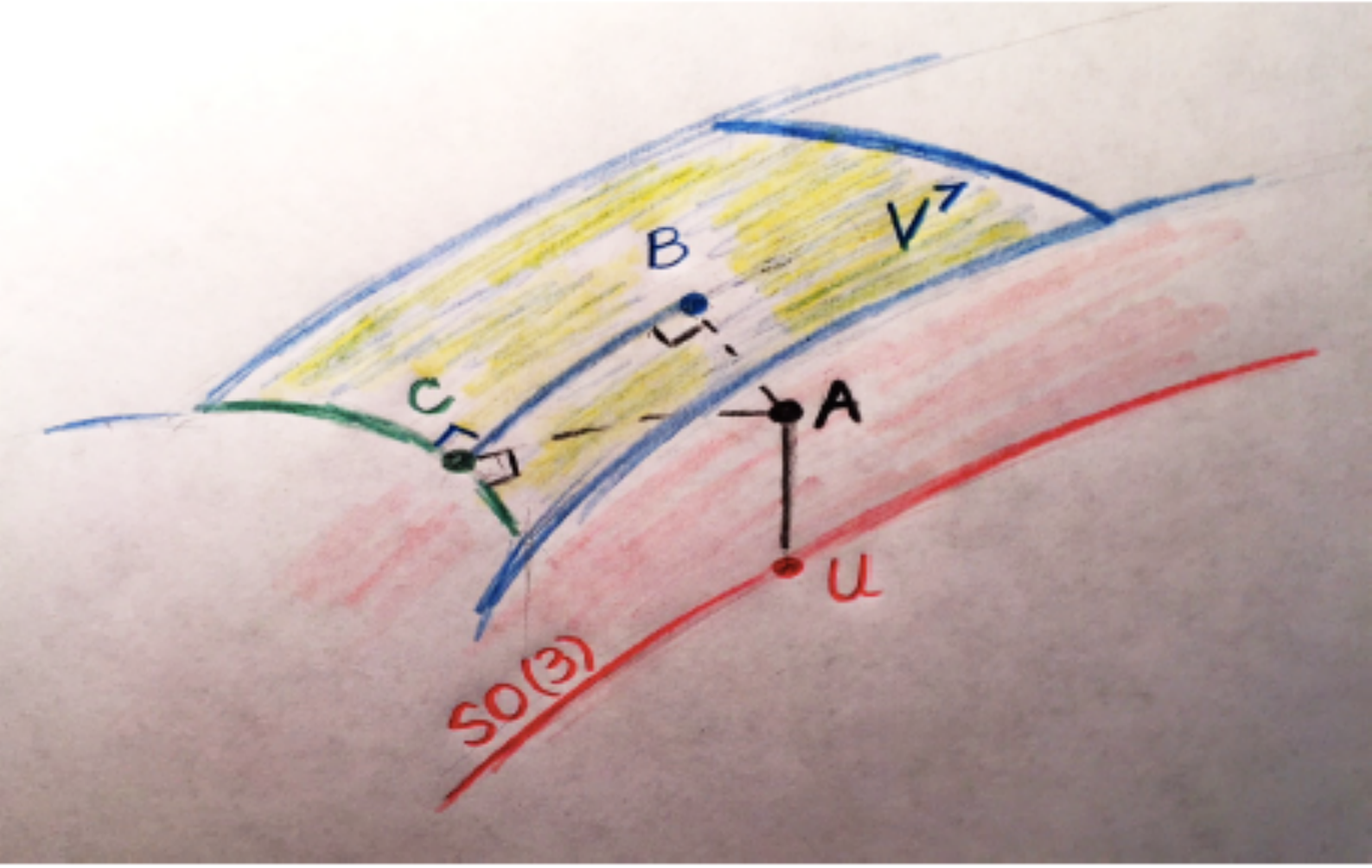}
\caption{{\textit{\textbf{\boldmath Nearest singular matrix}}}}
}
\end{figure}
In that figure, we start with a $3 \times 3$ matrix $A$ on $S^8(\sqrt{3})$, having positive
determinant and thus lying within the tubular neighborhood $N$ of $SO(3)$, with
$U$ its nearest orthogonal neighbor.
If $B$ is the nearest singular matrix on $S^8(\sqrt{3})$ to $A$, and $C$ is the nearest rank 1
matrix on $S^8(\sqrt{3})$ to $B$, then $C$ will also be the nearest rank 1 matrix there to $A$.

\bigskip

\medskip

\centerline{\Large{\textbf{\boldmath{Principal component analysis}}}}

\medskip

\vspace{-10pt}

Consider the singular value decomposition
$A = W D V^{-1}$
of an $n \times k$ matrix $A$,
where $V$ is a $k \times k$ orthogonal matrix, where $D$ is an $n \times k$ diagonal matrix,
$$D = \diag (d_1 , d_2 ,\ldots, d_r),\quad \mbox{with}\  d_1 \ge d_2 \ge\cdots\ge  d_r \ge 0,$$
with $r = \min(k, n)$, and where $W$ is an $n \times n$ orthogonal matrix.

Suppose that the rank of $A$ is $s\le r = \min(k, n)$, and that $s' < s$.
Then from the Eckart-Young theorem, we know that the nearest $n\times k$ matrix $A'$ of
rank $\le s' < s$ is given by $A' = W D' V^{-1}$, with $W$ and $V$ as above, and with
$$D' = \diag(d_1 , d_2 ,\ldots, d_{s'}).$$
The image of $A'$ has the orthonormal basis $\{\w_1 , \w_2 ,\ldots, \w_{s'}\}$, which are the first
$s'$ columns of the matrix $W$.

The columns of $W$ are the vectors $\w_1 ,\w_2 ,\ldots$, and are known as the
\textit{\textbf{principal components}} of the matrix $A$ , and the first $s'$ of them span the image
of the best rank $s'$ approximation to $A$.

If the matrix $A$ is used to collect a family of data points, and these data points
are listed as the columns of $A$, then the orthonormal columns of $W$ are regarded
as the principal components of this family of data points.

But if the data points are listed as the rows of $A$, then it is the orthonormal
columns of $V$ which serve as the principal components.

\noindent\textbf{Remark}.\ \ Principal Component Analysis began with Karl Pearson [1901]. He wanted to find
the line or plane of closest fit to a system of points in space, in which the measurement
of the locations of the points are subject to errors in any direction. 

\begin{wrapfigure}[11]{r}{0.52\textwidth}
\vspace{-15pt}
\center{\includegraphics[width=0.50\textwidth]{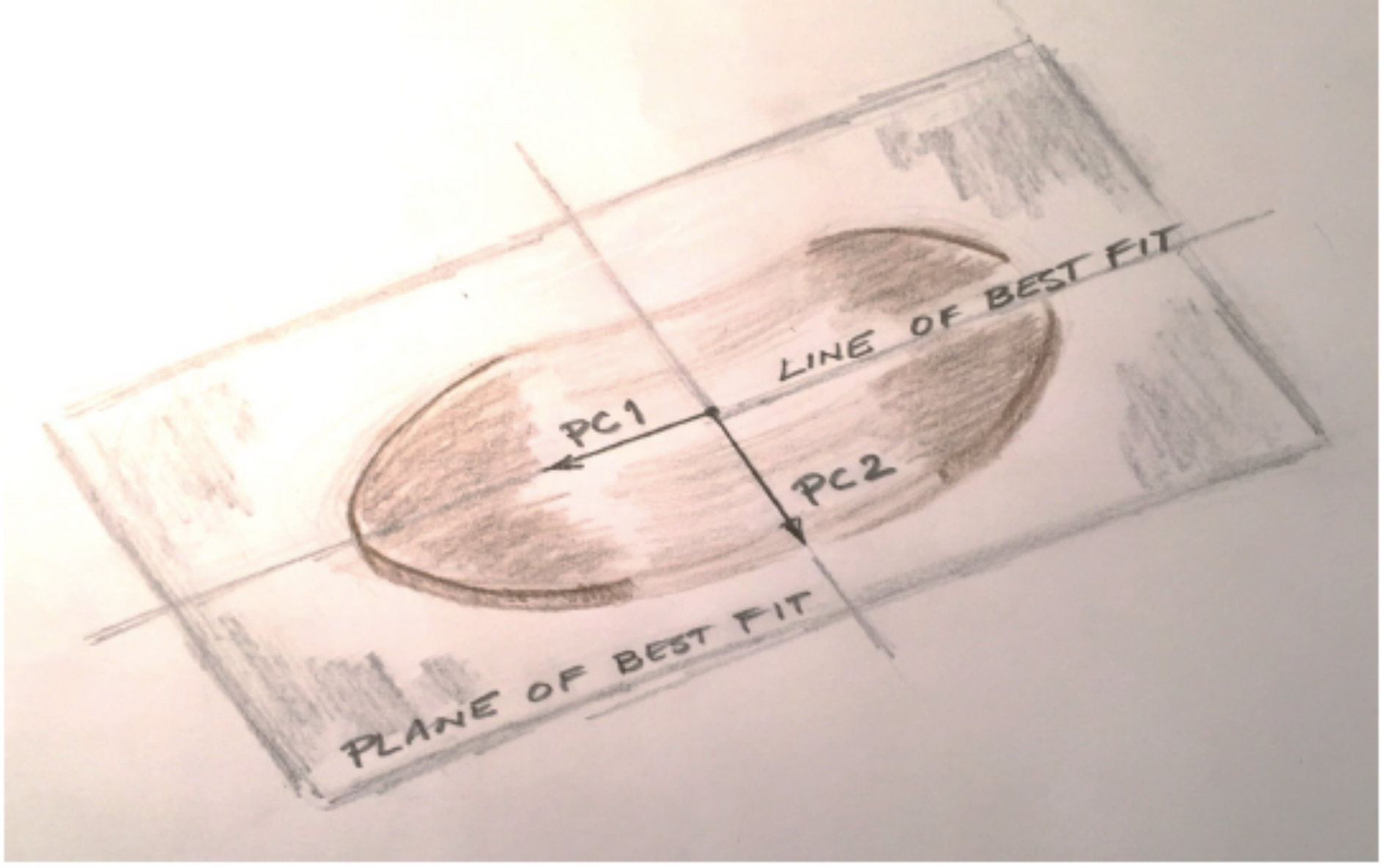}
\caption{{\textit{\textbf{\boldmath Principal components 1 and 2}}}}
}
\end{wrapfigure}
His key observation
was that to achieve this, one should seek to minimize the sum of the squares of the
perpendicular distances from all the points to the proposed line or plane of best fit.
The best fitting line is what we now view as the first principal component, described earlier (Figure 17). 

The actual term ``\textit{principal component}'' was introduced by Harold Hotelling [1933].

For further reading about the history of these matrix decompositions, we recommend
Horn and Johnson [1991], pages 134--140, and Stewart [1993] as excellent resources.

\vfill
\eject

\bigskip

\centerline{\LARGE{\textbf{\textit{Applications of nearest orthogonal neighbor}}}}

\bigskip

\centerline{\Large{\textbf{\boldmath{The orthogonal Procrustes problem}}}}

\vspace{-5pt}

Let $P = \{\p_1 , \p_2 , \ldots , \p_k\}$ and $Q = \{\q_1 , \q_2 , \ldots , \q_k\}$ be two ordered sets of
points in Euclidean $n$-space $\reals^n$. We seek a rigid motion $U$ of $n$-space which moves
$P$ as close as possible to $Q$, in the sense of minimizing the \textit{\textbf{disparity}} 
$d_1^2+d_2^2+\cdots+d_k^2$ between $U(P)$ and $Q$, where
$d_i = \Vert U(\p_i) - \q_i\Vert$.

It is easy to check that if we first translate the sets $P$ and $Q$ to put their 
centroids
at the origin, then this will guarantee that the desired rigid motion $U$ also fixes the
origin, and so lies in $O(n)$.
We assume this has been done, so that the sets $P$ and $Q$ have their centroids
at the origin.

Then we form the $n \times k$ matrices $A$ and $B$ whose columns are the vectors
$\p_1 , \p_2 , \ldots , \p_k$ and $\q_1 , \q_2 , \ldots, \q_k$, and we seek the matrix $U$ in 
$O(n)$ which
minimizes the disparity $d_1^2+d_2^2+\cdots+d_k^2 = \Vert U A - B \Vert^2$ between $U(P)$
and $Q$.

We start by expanding
$$\langle U A - B \,,\, U A - B \rangle = \langle U A \,,\, U A \rangle - 2 \langle U A \,,\, B \rangle + 
\langle B , B \rangle.$$
Now $\langle U A \,,\, U A \rangle = \langle A , A \rangle$ which is fixed, and likewise 
$\langle  B , B \rangle$ is fixed,
so we want to \textit{maximize} the inner product $\langle U A \,,\, B \rangle$ by appropriate choice of 
$U$ in $O(n)$. But
$$\langle U A \,,\, B \rangle = \langle U \,,\, B A^T \rangle,$$
and so, reversing the above steps, we want to \textit{minimize} the inner product
$$\langle U -B A^T \,,\, U - B A^T \rangle,$$
which means that we are seeking the orthogonal transformation $U$ which is
closest to $B A^T$ in the space of $n\times n$ matrices.

The above argument was given by Peter Sch\"onemann [1966] in his PhD thesis at
the University of North Carolina.

When $n \ge 3$, we don't have a simple explicit formula for $U$, but it is the
orthogonal factor in the polar decomposition
$$B A^T = U P = P' U.$$
Visually speaking, if we scale $B A^T$ to lie on the round $n^2 - 1$ sphere of radius
$\sqrt{n}$ in\break $n^2$-dimensional Euclidean space $\reals^{n^2}$, then $U$ is at the center of the
cross-sectional cell in the tubular neighborhood of $O(n)$ which contains $B A^T$,
and is unique if $\det (B A^T) \ne 0$.

\bigskip

\centerline{\Large{\textbf{\boldmath{A least squares estimate of satellite attitude}}}}

\medskip

\vspace{-10pt}

Let $P = \{\p_1 , \p_2 , \ldots , \p_k\}$ be unit vectors in 3-space which represent the direction
cosines of $k$ objects observed in an earthbound fixed frame of reference, and\break
$Q = \{\q_1 , \q_2 , \ldots , \q_k\}$ the direction cosines of the same $k$ objects as observed in a
satellite fixed frame of reference. Then the element $U$ in $SO(3)$ which minimizes
the disparity between $U(P)$ and $Q$ is a least squares estimate of the rotation matrix
which carries the known frame of reference into the satellite fixed frame at any
given time. See Wahba [1966].

Errors incurred in computation of $U$ can result in a loss of orthogonality, and be
compensated for by moving the computed $U$ to its nearest orthogonal neighbor.

\bigskip

\medskip

\centerline{\Large{\textbf{\boldmath Procrustes best fit of anatomical objects}}}

\medskip

\vspace{-10pt}

The challenge is to compare two similar anatomical objects: two skulls, two teeth,
two brains, two kidneys, and so forth.

Anatomically corresponding points
(\textit{landmarks}) are chosen on the two objects, say the ordered set of points $P = \{\p_1 ,\p_2 , \ldots , \p_k\}$
on the first object, and the ordered set of points $Q = \{\q_1 , \q_2 ,\ldots, \q_k\}$ on the second
object.  They are translated to place their centroids at the origin,
and then the Procrustes procedure is applied by seeking a rigid motion $U$ of 3-space so
as to minimize the disparity 
$d_1^2 + d_2^2 + \cdots + d_k^2$ between $U(P)$ and $Q$ , where
$d_i = \Vert U(\p_i) - \q_i \Vert$.
\begin{figure}[h!]
\center{\includegraphics[height=144pt]{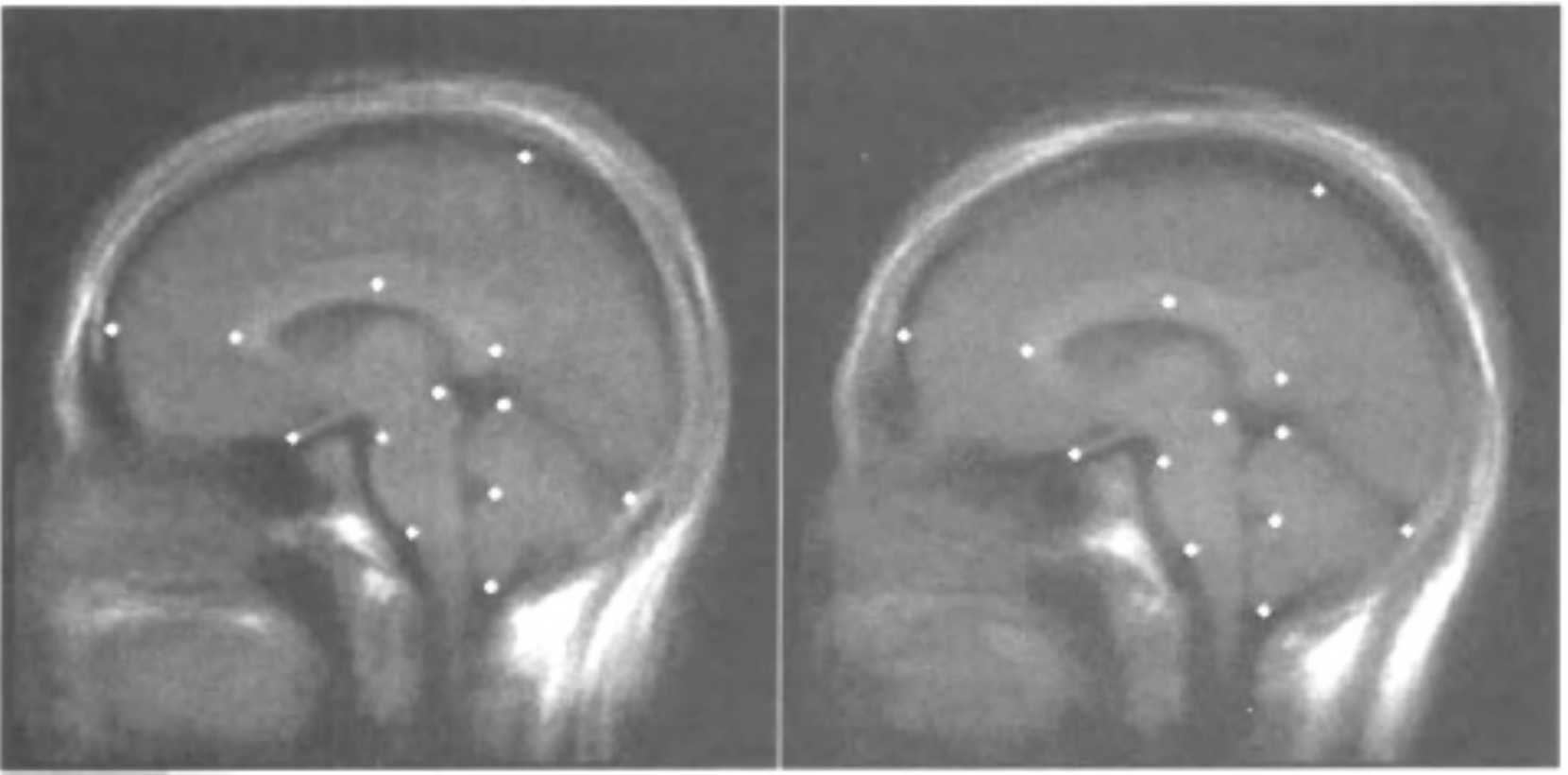}
\caption{{\textit{\textbf{\boldmath Brain scans}}}}
}
\end{figure}

In Figure 18, the left brain slice is actually an average over a group of doctors,
and the right slice an average over a group of patients, each with 13 corresponding
landmark points, from the paper by Bookstein [1997].

If size is not important in the comparison of two shapes, then it can be factored out by
scaling the two sets of landmarks, $P = \{\p_1 ,\p_2 , \ldots , \p_k\}$ and $Q = \{\q_1 , \q_2 ,\ldots, \q_k\}$,
so that 
$\Vert\p_1\Vert^2 + \cdots + \Vert\p_k\Vert^2 = \Vert\q_1\Vert^2 +\cdots + \Vert\q_k\Vert^2$.

For modifications which allow comparison of any number of shapes at the same
time, see for example Rohlf and Slice [1990].

The effectiveness of this Procrustes comparison naturally depends on appropriate
choice and placement of the landmark points, and leads one to seek an alternative
approach which does not depend on this. To that end, see Lipman, Al-Aifari and Daubechies [2013]
in which the authors propose a continuous Procrustes distance, and then prove that it
provides a metric for the space of ``shapes'' of two-dimensional surfaces embedded
in three-space.

\bigskip

\medskip

\centerline{\Large{\textbf{\boldmath{Facial recognition and eigenfaces}}}}

\medskip

\vspace{-10pt}

We follow Sirovich and Kirby [1987]
in which the principal components of the data base matrix of facial pictures are
suggestively called \textit{\textbf{eigenpictures}}.

The authors and their team assembled a file of 115 pictures of undergraduate students
at Brown University. Aiming for a relatively homogeneous population, these students
were all smooth-skinned caucasian males. The faces were lined up so that the same
vertical line passsed through the symmetry line of each face, and the same horizontal
line through the pupils of the eyes. Size was normalized so that facial width was the
same for all images.

Each picture contained $128 \times 128 = 2^{14} = 16,384$ pixels, with a grey scale
determined at each pixel. So each picture was regarded as a single vector
$\varphi^{(n)}$,\break  $n = 1, 2,\ldots, 115$, called a \textit{\textbf{face}}, in a vector space of dimension $2^{14}$.

\begin{wrapfigure}[10]{r}{0.27\textwidth}
\vspace{-30pt}
\center{\includegraphics[width=0.23\textwidth]{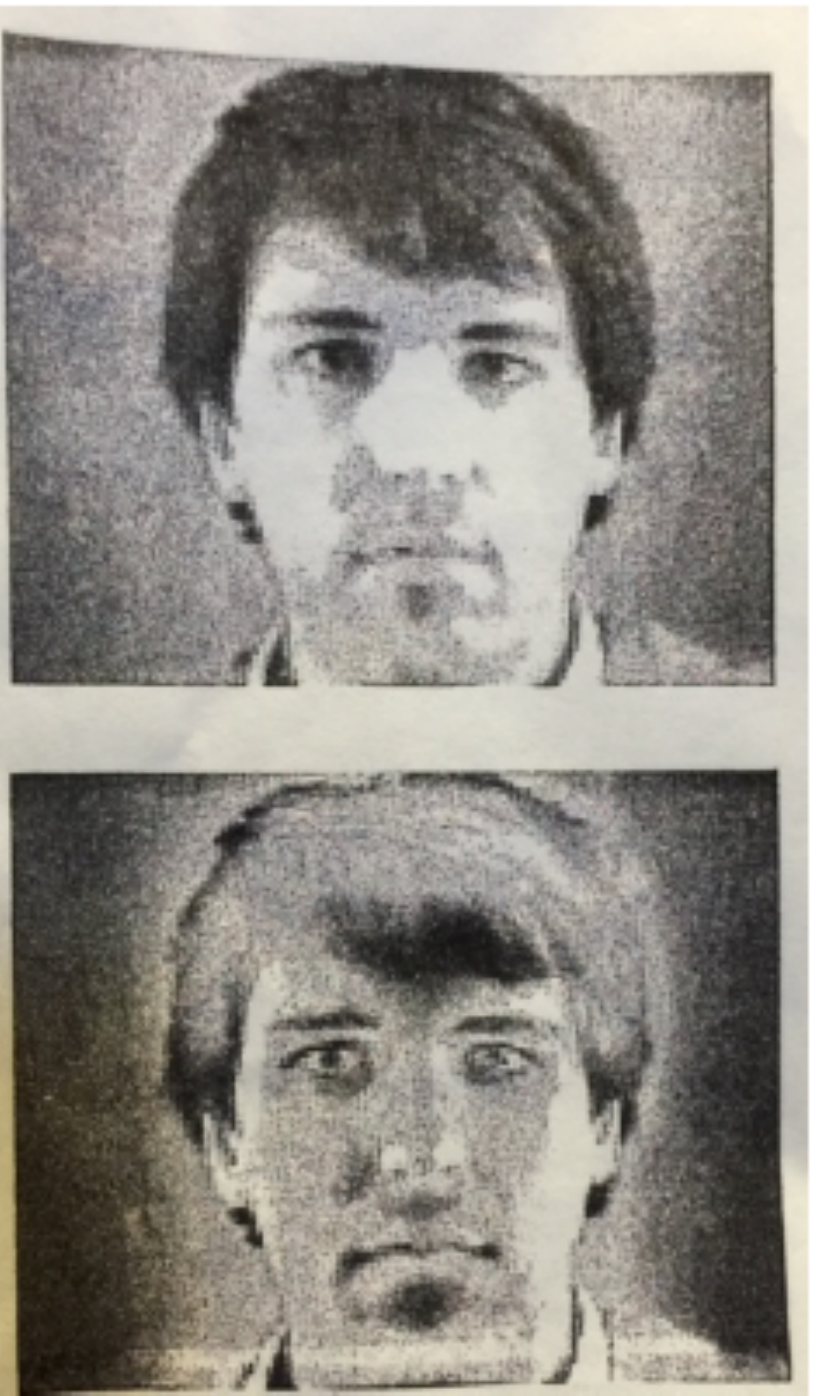}
\caption{{\textit{\textbf{\boldmath Sample face and caricature}}}}
}
\end{wrapfigure}
The challenge was to find a low-dimensional subspace of best fit to these 115 faces,
so that a person could be sensibly recognized by the projection of his picture into this
subspace.

To make sure that the subspace passes through the origin (i.e., is a linear rather than
affine subspace), the data is adjusted so that its average is zero, as follows.

Let $\displaystyle\langle\varphi\rangle = (1/M) \sum_{n=1}^M \varphi^{(n)}$
 be the average face, where $M = 115$, and then let
$\phi^{(n)} = \varphi^{(n)} - \langle\varphi\rangle$
be the deviation of each face from the average. The authors refer to each such
deviation $\phi$ as a \textit{\textbf{caricature}}.
Figure 19 shows a sample face, and its caricature.

The collection of caricatures $\phi^{(n)},\  n = 1, 2, \ldots, 115$ was then regarded as a $2^{14} \times 115$
matrix $A$, with each caricature appearing as a column of $A$.

If the singular value decomposition of $A$ is
$A = W D V^{-1}$,
with $W$ a $2^{14}\times 2^{14}$ orthogonal matrix,
$$D = \diag(d_1 , d_2 ,\ldots, d_{115})\quad \mbox{with}\ d_1\ge d_2 \ge\cdots \ge d_{115} \ge 0$$
a $2^{14}\times  115$ diagonal matrix, and $V$ a $115 \times 115$ orthogonal matrix,
then
the orthonormal columns $\w_1 , \w_2 ,\ldots, \w_{2^{14}}$ of $W$ are the principal components of
the matrix $A$.

\begin{wrapfigure}[16]{r}{0.47\textwidth}
\vspace{-30pt}
\center{\includegraphics[width=0.43\textwidth]{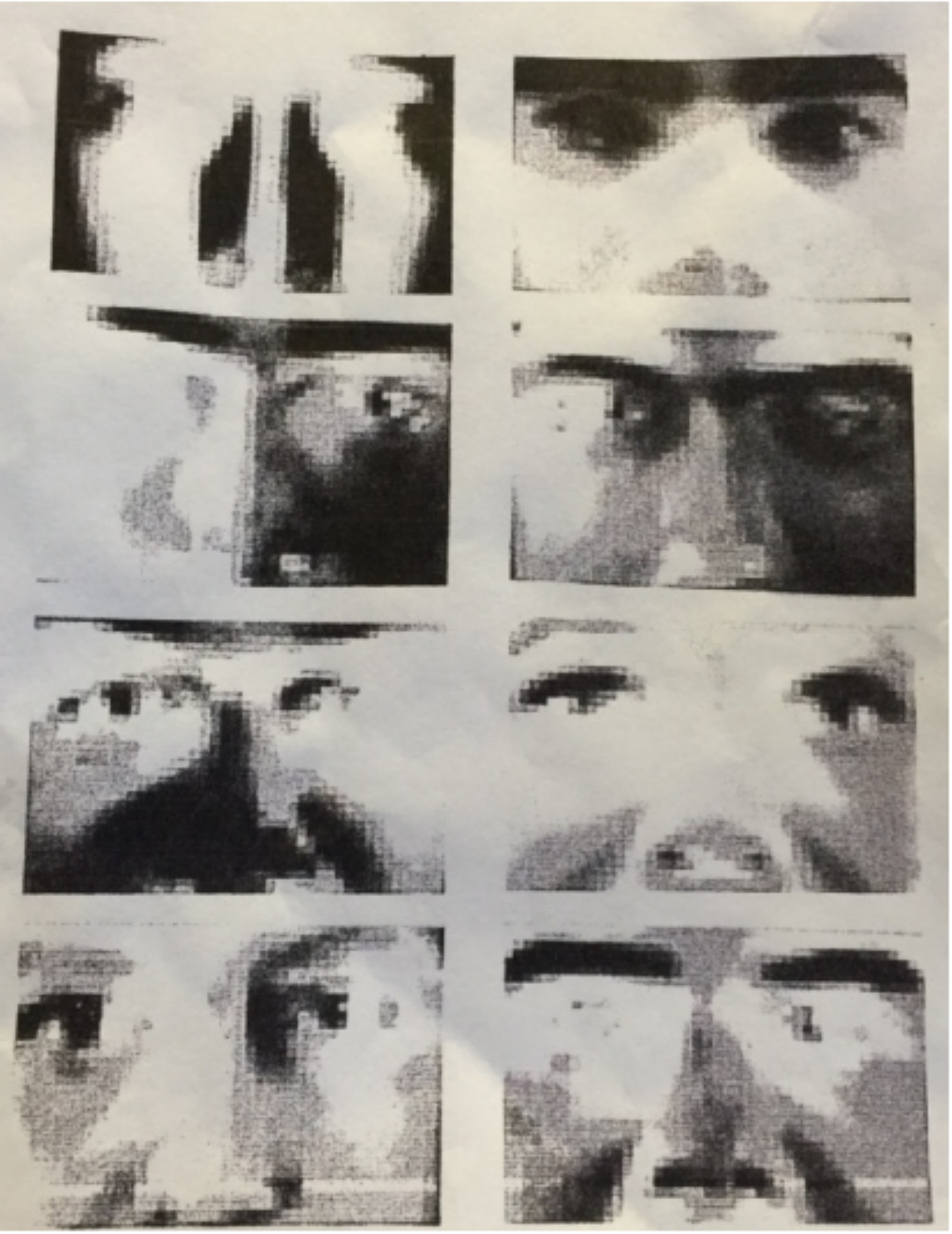}
\caption{{\textit{\textbf{\boldmath First eight eigenfaces}}}}
}
\end{wrapfigure}
It was found that the first 100 principal components of $A$ span a subspace
sufficiently large to recognize any of the faces $\varphi^{(n)}$ by projecting its caricature
into this subspace and then adding back the average face:
$$\varphi^{(n)} \sim  \langle\varphi\rangle+\sum_{k=1}^{100}  \langle \phi^{(n)}\,,\,\w_k\rangle\,\w_k.$$

Figure 20 shows the first eight eigenpictures starting at the upper left, moving to the
right, and ending at the lower right, in which each picture is cropped to focus
on the eyes and nose. Since the eigenpictures can have negative entries, a
constant was added to all the entries to make them positive for the purpose of viewing

Figure 21 shows a sample face, correspondingly cropped,
\begin{figure}[h!]
\center{\includegraphics[height=120pt]{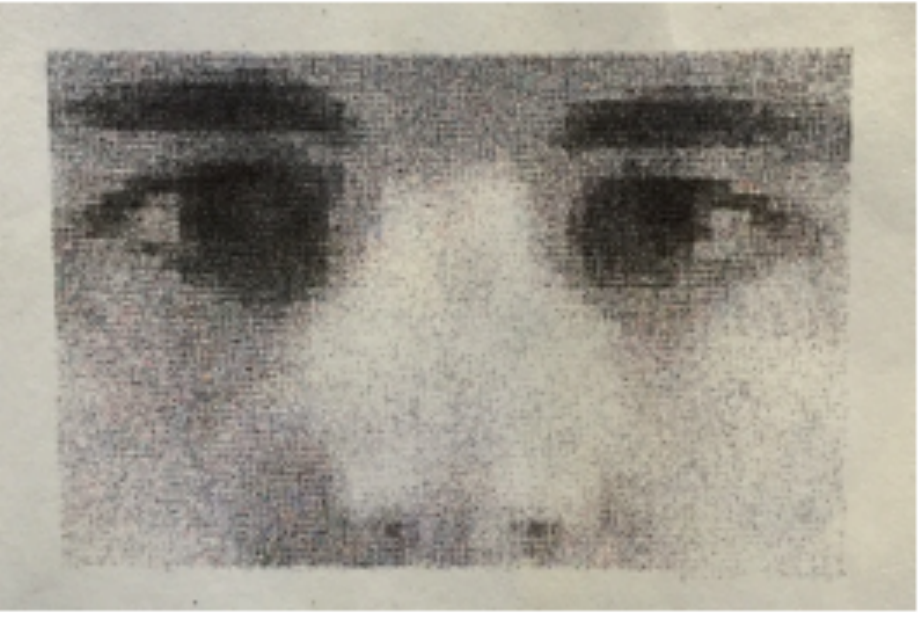}
\caption{{\textit{\textbf{\boldmath Cropped sample face}}}}
}
\end{figure}
and Figure 22 shows the approximations to that sample face, using 10, 20, 30 and 40
eigenpictures.
\begin{figure}[h!]
\center{\includegraphics[height=180pt]{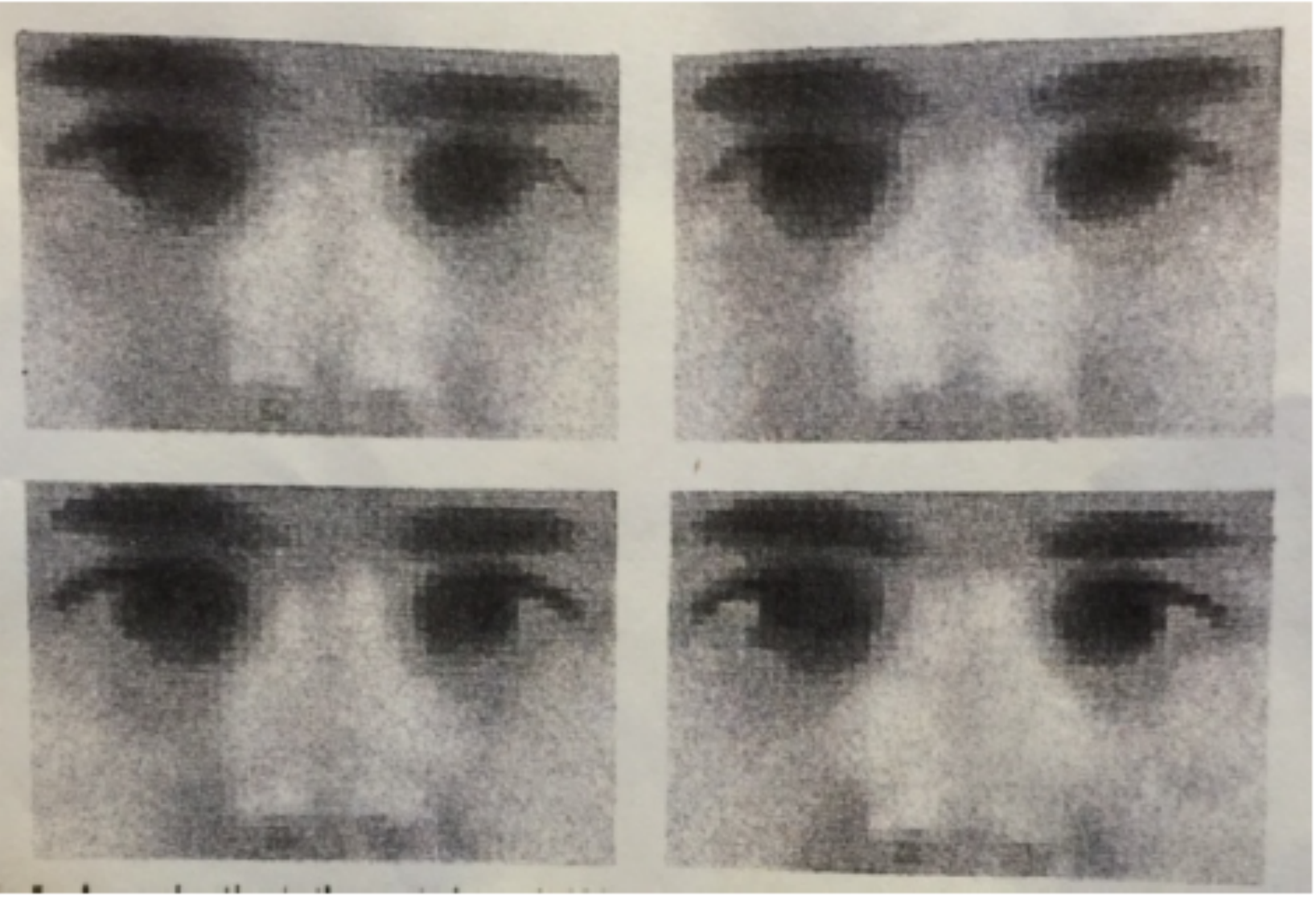}
\caption{{\textit{\textbf{\boldmath Approximations to sample face}}}}
}
\end{figure}

After working with the initial group of 115 male students, the authors tried out
the recognition procedure on one more male student and two females, using
40 eigenpictures, with errors of 7.8\%, 3.9\%, and 2.4\% in these three cases.

\noindent\textbf{\textit{Remarks.}}
\vspace{-10pt}
\begin{enumerate}
\item[(1)] In the pattern recognition literature, the Principal Component Analysis method used in this paper is
also known as the Karhunen-Loeve expansion.
\item[(2)] Another very informative and nicely written paper on this approach to facial
recognition is Turk and Pentland [1991]. The section of this paper on \textit{Background
and Related Work} is a brief but very interesting survey of alternative approaches to
computer recognition of faces.
An overview of the literature on face recognition is given in
Zhao et al [2003].
\end{enumerate}

\bigskip

\centerline{\Large{\textbf{\boldmath{Principal component analysis applied to}}}}
\vspace{3pt}
\centerline{\Large{\textbf{\boldmath{ interest rate term structure}}}}

\medskip

\vspace{-10pt}

How does the interest rate of a bond vary with respect
to its \textit{term}, meaning time to maturity? The answer involves one of the oldest and best known
applications of Principal Components Analysis (PCA) to the field of economics
and finance, originating in the work of Litterman and Scheinkman [1991].

To begin, economists plot the interest rate for a given bond against a variety of
different maturities, and call this a \textit{\textbf{yield curve}}. Figure 23 shows such a curve for US
Treasury bonds from an earlier date, when interest rates were higher than they
are now.

\begin{figure}[h!]
\center{\includegraphics[height=180pt]{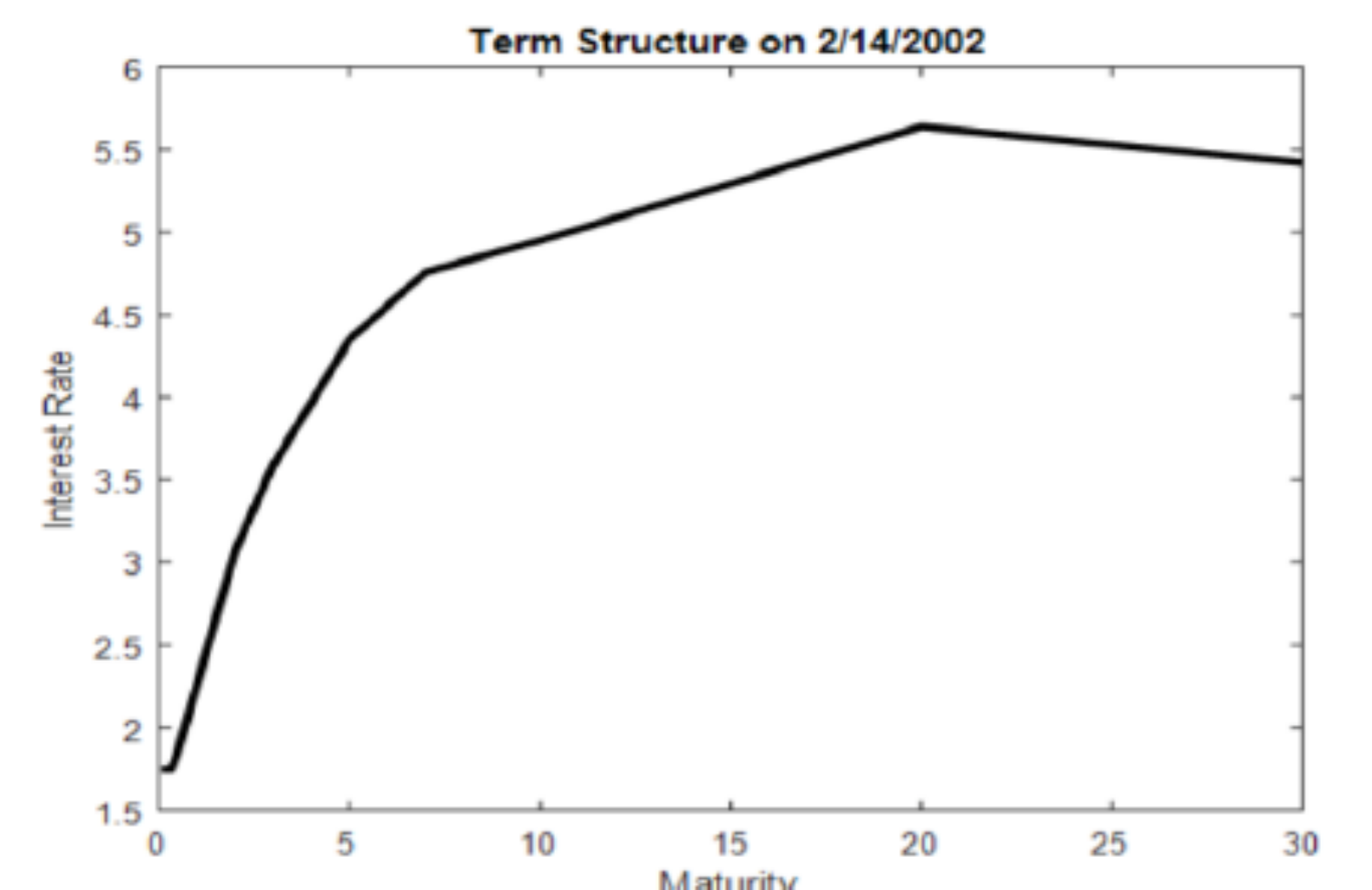}
\caption{{\textit{\textbf{\boldmath Yield curve}}}}
}
\end{figure}
\vspace{-8pt}
Predicting the relation shown by such a curve can be crucial for investors
trying to determine which assets to invest in, and for governments who wish to
determine the best mix of Treasury maturities to auction on any given day.
For this reason, a number of investigators have tried to understand whether there
are common factors embedded in the term structure. In particular, identifying
whether there are factors which affect all interest rates equally, or which affect
interest rates for bonds of certain maturities but not of others, is important for
understanding how the term structure behaves.

To help understand how these questions are answered, we replicated the methodology in the Litterman and
Scheinkman paper, using a newer data set which gives the daily interest rate term
structure for US Treasury bonds over a long span of time, 2,751 days between 2001 and 2016. 
For each of these days, we recorded the interest rates for bonds of
11 different maturities: 1, 3 and 6 months, and 1, 2, 3, 5, 7, 10, 20 and 30 years.
Each data vector is an 11-tuple of interest rates, which we collect as the rows of a
$2,751 \times 11$ matrix.

The average of the rows is  depicted graphically in Figure 24.
\begin{figure}[h!]
\center{\includegraphics[height=164pt]{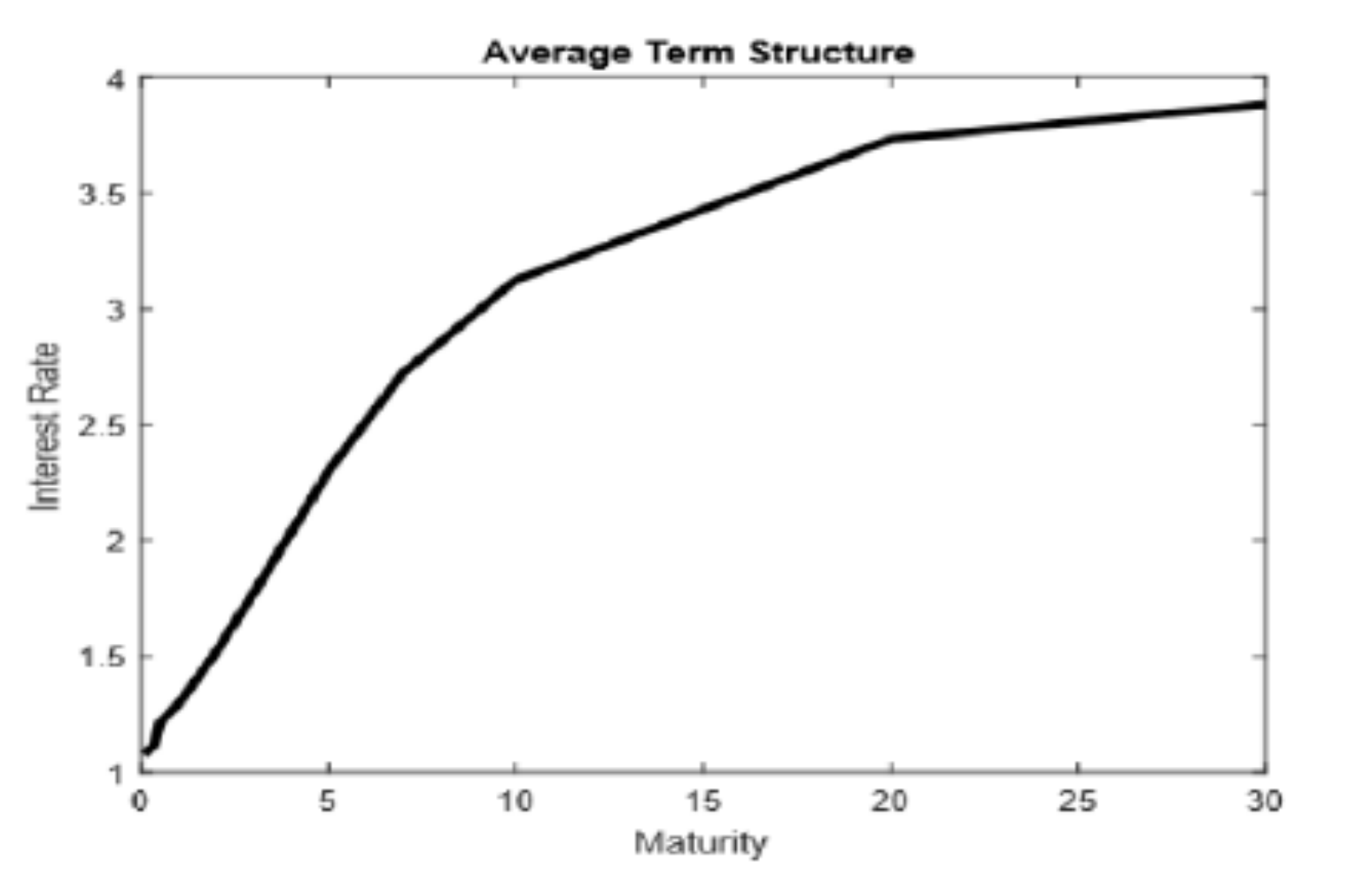}
\caption{{\textit{\textbf{\boldmath Average yield curve}}}}
}
\end{figure}
We subtracted this average from each of the rows, and called the
resulting matrix $A$. The rows of $A$ are our \textit{adjusted data vectors}, which now
add up to zero.

Let $A = W D V^{-1}$ be the singular value decomposition of $A$, where $V$ is an $11 \times 11$
orthogonal matrix, $D$ is a $2,751 \times 11$ diagonal matrix, and $W$ is a $2,751 \times 2,751$
orthogonal matrix. Since the data points are the \textit{rows} of $A$, the principal components
are the 11 orthonormal columns of $V$.

These principal components reveal the line of best fit, the plane of best fit, the
3-space of best fit, and
so forth for our 2,751 data points. They were obtained using the PCA package of MATLAB.
The first three principal components are shown graphically in Figure 25.
\begin{figure}[h!]
\center{\includegraphics[height=164pt]{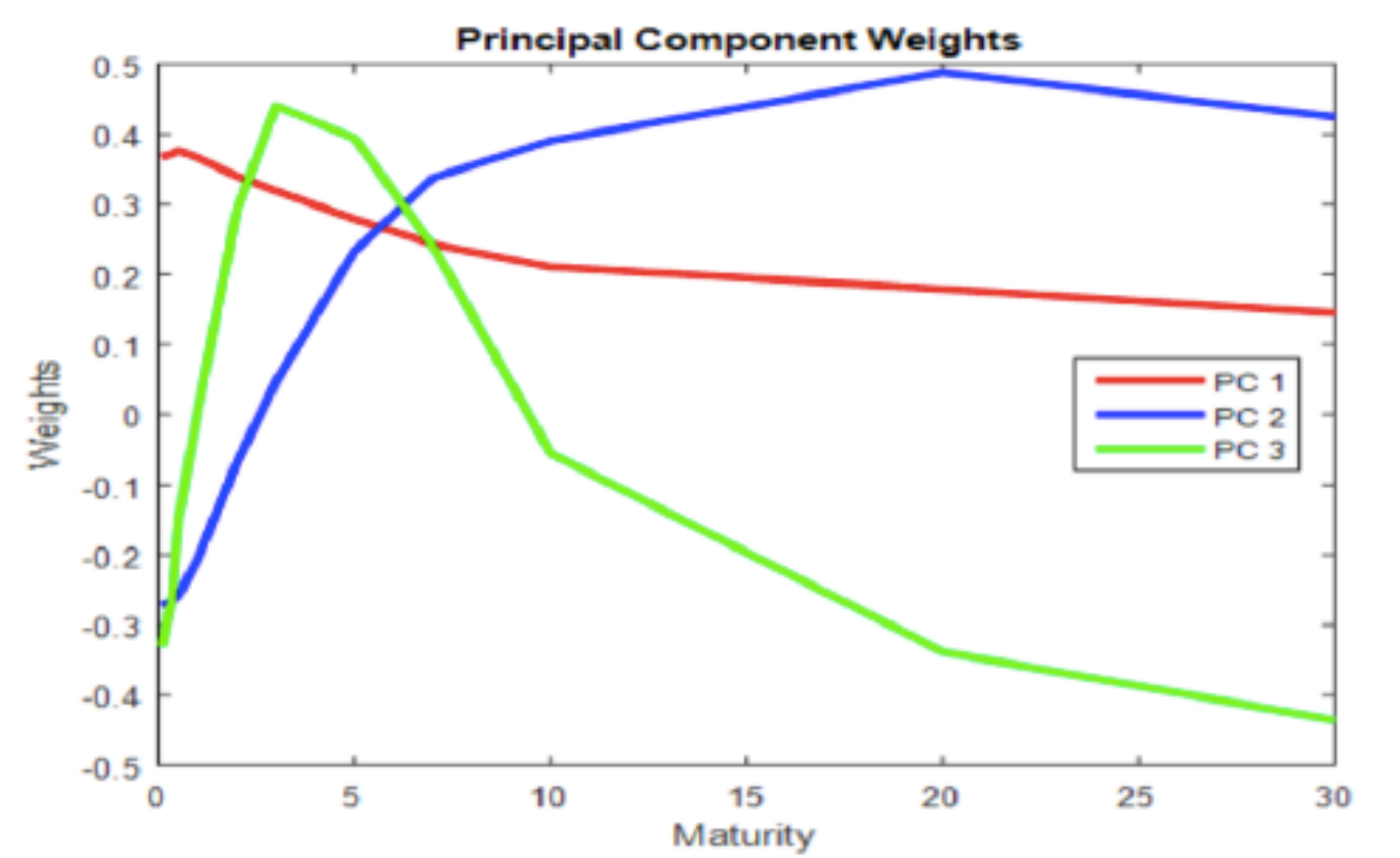}
\caption{{\textit{\textbf{\boldmath Principal components of $A$}}}}
}
\end{figure}
The first principal component is more constant than the other two, and captures
the fact that most of the variation in term structures comes from changes which
affect the levels of all yields.

The second most important source of variation in term structure comes from the
second principal component, which reflects changes that most affect yields on
bonds of longer maturities, while the third principal component reflects changes
that affect medium term yields the most.
These features of the first three principal components were called \textit{level}, \textit{steepness},
and \textit{curvature} in the foundational paper by Litterman and Scheinkman.

In Figure 26, the black curve is the term structure on 2/14/2002, duplicating the
first figure in this section. We subtract the average term structure from this
particular one, project the difference onto the one-dimensional subspaces spanned
in turn by the first three principal components, and show these projections below in
red, blue and green. Finally, we sum up these three projections, add back the average
term structure, show the result in purple, and see how closely this purple curve
approximates the black curve we started with.
\begin{figure}[h!]
\center{\includegraphics[height=180pt]{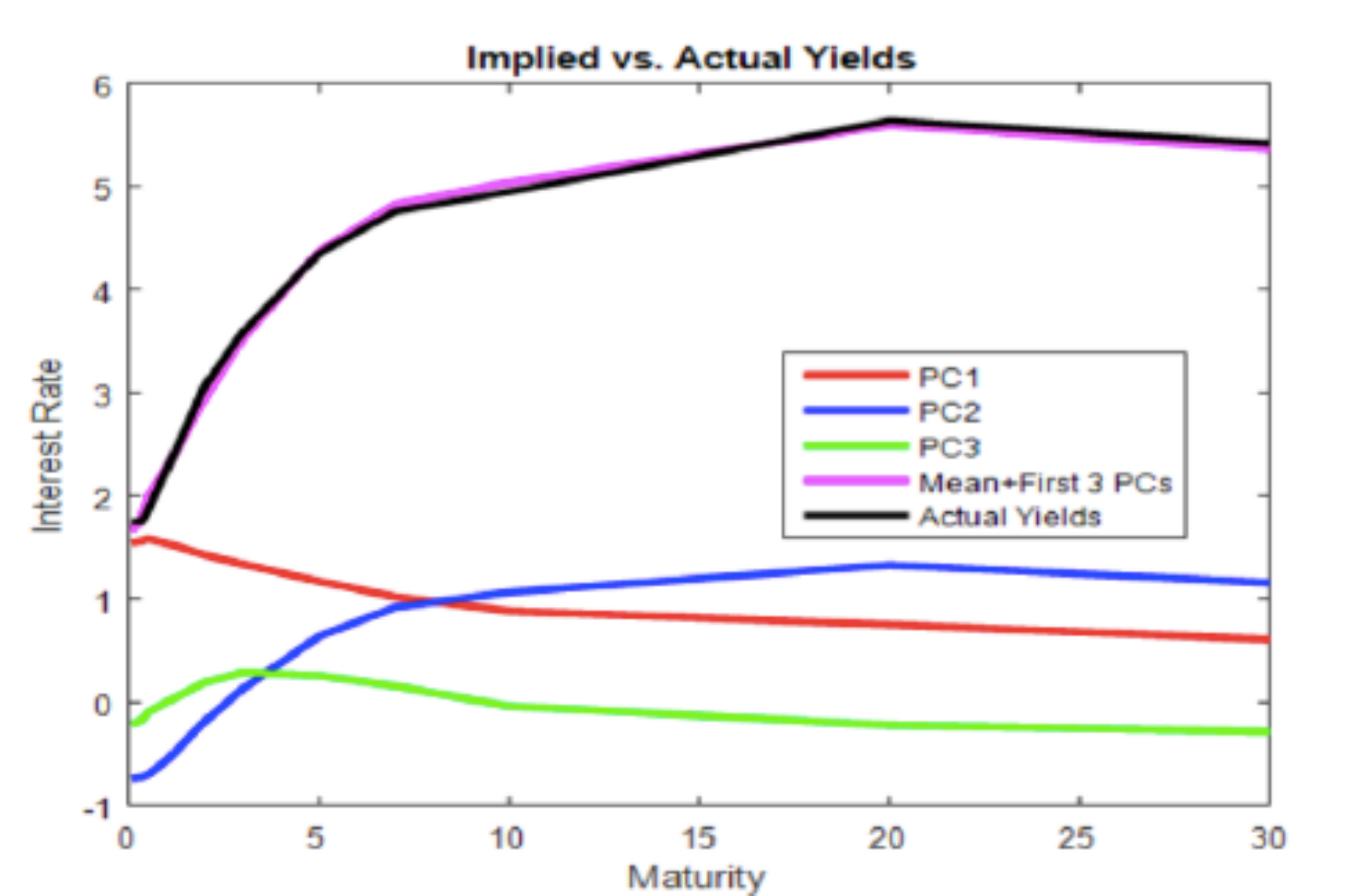}
\caption{{\textit{\textbf{\boldmath Approximation of a yield curve by its first three principal components}}}}
}
\end{figure}

\bigskip

\centerline{\bf REFERENCES\rm}

\addtolength{\baselineskip}{-1pt}

\vspace{-10pt}

\begin{description}
\item[1873] E.\ Beltrami, \textit{Sulle funzioni bilineari}, Giornali di Mat.\ ad Uso degli Studenti
Delle Universita, 11, 98--106
\item[1874a] C.\ Jordan, \textit{Memoire sur les formes bilineares}, J.\ Math.\ Pures Appl.,
2nd series, 19, 35--54
\item[1874b] C.\ Jordan, \textit{Sur la reduction des formes bilineares}, Comptes Rendus de
l'Academie Sciences, Paris 78, 614--617
\item[1901] Karl Pearson, \textit{On Lines and Planes of Closest Fit to Systems of Points
in Space}, Philosophical Magazine 2, 559--572.
\item[1902] L.\ Autonne, \textit{Sur les groupes lineaires, reels et orthogonaux}, Bull.\ Soc.\ Math.\ France, 30, 121--134
\item[1907] E.\ Schmidt, \textit{Zur Theorie des linearen und nichlinearen Integral gleichungen,
I Teil. Entwicklung willkurlichen Funktionen nach System vorgeschriebener},
Math.\ Ann.\ 63, 433 -- 476
\item[1933] H.\ Hotelling, \textit{Analysis of a complex of statistical variables into principal components},
J.\  Ed.\ Psych, 24, 417 -- 441 and 498 -- 520
\item[1936] C.\ Eckart and G.\ Young, \textit{The approximation of one matrix by another of lower
rank}, Psychometrika, I, 211 -- 218.
\item[1947] Kari Karhunen, \textit{Uber lineare Methoden in der Wahrscheinlichkeitsrechnung},
Ann.\ Acad.\ Sci.\ Fennicae, Ser.\ A.\ I.\ Math-Phys 37, 1--79.
\item[1961] John Francis, \textit{The QR transformation}, parts I and II, Computer J.\ Vol.\ 4,
265--272 and 332--345.
\item[1962] Vera Kublanovskaya, \textit{On some algorithms for the solution of the complete
eigenvalue problem}, USSR Comput.\ Math.\ and Math.\ Physics. vol 1, 637--657.
\item[1966] Peter H.\ Sch\"onemann, \textit{A generalized solution of the orthogonal Procrustes
problem}, Psychometrika, Vol.\ 31, No.\ 1, March, 1 -- 10.
\item[1966] Grace Wahba, \textit{A Least Squares Estimate of Satellite Attitude}, SIAM Review,
Vol.\ 8, No.\ 3,  384 -- 386.
\item[1976] Wolfgang Kabsch, \textit{A solution for the best rotation to relate two sets of
vectors}, Acta Crystallographica 32, 922, with a correction in 1978,
\textit{A discussion of the solution for the best rotation to relate two sets of
vectors}, Ibid, A-34, 827--828.
\item[1980] H.\ Blaine Lawson, \textit{Lectures on Minimal Submanifolds}, Publish or Perish Press.
\item[1985] Roger Horn and Charles Johnson, {\it Matrix Analysis}, Cambridge University Press.
\item[1986] Nicholas J.\ Higham, \textit{Computing the polar decomposition -- with applications},
SIAM J.\ Sci.\ Stat.\ Comput.\ Vol.\ 7, No.\ 4 October, 1160 -- 1174.
\item[1987] L.\ Sirovich and M.\ Kirby, \textit{Low-dimensional procedure for the characterization
of human faces}, J.\ Optical Society of America, Vol.\ 4, No.\ 3, 519 -- 524,
\item[1990] F.\ James Rohlf and Dennis Slice, \textit{Extensions of the Procrustes method for the
optimal superimposition of landmarks}, Syst.\ Zool.\ 39 (1), 40 -- 59.
\item[1991] Roger Horn and Charles Johnson, \textit{Topics in Matrix Analysis},
Cambridge University Press.
\item[1991] Robert Litterman and Jos\'e Scheinkman, \textit{Common factors affecting bond returns},
J.\ Fixed Income, June, 54 -- 61.
\item[1991] Matthew Turk and Alex Pentland, \textit{Eigenfaces for Recognition}, Journal of Cognitive
Neuroscience, Vol.\ 3, No.\ 1, 71 -- 86.
\item[1993] G.W.\ Stewart, \textit{On the early history of the singular value decomposition},
SIAM Review, Vol.\ 35, No.\ 4, 551 -- 566
\item[1996] Gene Golub and Charles Van Loan, \textit{Matrix Computations}, Third Edition,
Johns Hopkins University Press.
\item[1997] Fred L.\  Bookstein, \textit{Biometrics and brain maps: the promise of the Morphometric
Synthesis}, in S.\  Kowlow and M.\ Huerta, eds., \textit{Neuroinformatics: An Overview of
the Human Brain Project}, Progress in Neuroinformatics, Vol.\ 1, 203 -- 254.
\item[2000] Barry Cipra, \textit{The Best of the 20th Century: Editors Name Top 10 Algorithms},
SIAM News, Vol.\ 33, No.\ 4, 1--2.
\item[2002] Allen Hatcher, \textit{Algebraic Topology}, Cambridge University Press.
\item[2003] W.\ Zhao, R.\ Chellappa, P.J.\ Phillips and A.\ Rosenfeld, \textit{Face Recognition:
A Literature Survey}, ACM Computing Surveys, Vol.\ 35, No.\ 4, 399 -- 458.
\item[2009] G.\ H.\ Golub and F.\ Uhlig, \textit{The QR algorithm: 50 years later its genesis by
John Francis and Vera Kublanovskaya and subsequent developments},
IMA J.\ Numer.\ Anal.\ Vol.\ 29, 467--485.
\item[2011] David S.\ Watkins, \textit{Francis's Algorithm}, American Mathematical Monthly
Vol.\ 118, May, 387--403.
\item[2013] Yaron Lipman, Reema Al-Aifari and Ingrid Daubechies, \textit{The continuous Procrustes
distance between two surfaces}, Comm.\ Pure Appl.\ Math. 66, 934 -- 964,
\end{description}

\noindent University of Pennsylvania\\
Philadelphia, PA  19104

\noindent
Dennis DeTurck: \it deturck@math.upenn.edu\rm\\
Amora Elsaify: \it aelsaify@wharton.upenn.edu\rm\\
Herman Gluck: \it gluck@math.upenn.edu\rm\\
Benjamin Grossmann: \it bwg25@drexel.edu\rm\\
Joseph Hoisington: \it jhois@math.upenn.edu\rm\\
Anusha M.\,Krishnan: \it anushakr@math.upenn.edu\rm\\
Jianru Zhang: \it jianruzh@math.upenn.edu\rm\\

\end{document}